\documentclass[11pt]{article}
\usepackage{newlfont,amsfonts,amssymb, oldgerm,amsmath,amsthm,amsgen,amscd,
}
\usepackage{epic}
\usepackage{eepic}
\usepackage{rotating}
\usepackage[dvips]{epsfig}
\begin{document}

\newcommand{\Q}{\ensuremath{\mathbb{H}}}
\newcommand{\N}{\ensuremath{\mathbb{N}}}
\newcommand{\Z}{\ensuremath{\mathbb{Z}}}
\newcommand{\C}{\ensuremath{\mathbb{C}}}
\newcommand{\K}{\ensuremath{\mathbb{K}}}
\renewcommand{\O}{\ensuremath{\mathcal{O}}}
\newcommand{\R}{\ensuremath{\mathbb{R}}}

\newcommand{\bcase}{\begin{case}}
\newcommand{\ecase}{\end{case}}
\newcommand{\setcase}{\setcounter{case}{0}}
\newcommand{\bclaim}{\begin{claim}}
\newcommand{\eclaim}{\end{claim}}
\newcommand{\setclaim}{\setcounter{claim}{0}}
\newcommand{\bstep}{\begin{step}}
\newcommand{\estep}{\end{step}}
\newcommand{\setstep}{\setcounter{step}{0}}
\newcommand{\bhlem}{\begin{hlem}}
\newcommand{\ehlem}{\end{hlem}}
\newcommand{\sethlem}{\setcounter{hlem}{0}}

\newcommand{\bleer}{\begin{leer}}
\newcommand{\eleer}{\end{leer}}
\newcommand{\bde}{\begin{de}}
\newcommand{\ede}{\end{de}}
\newcommand{\ul}{\underline}
\newcommand{\ol}{\overline}
\newcommand{\tbf}{\textbf}
\newcommand{\mc}{\mathcal}
\newcommand{\mb}{\mathbb}
\newcommand{\mf}{\mathfrak}
\newcommand{\bs}{\begin{satz}}
\newcommand{\es}{\end{satz}}
\newcommand{\btheo}{\begin{theo}}
\newcommand{\etheo}{\end{theo}}
\newcommand{\bfolg}{\begin{folg}}
\newcommand{\efolg}{\end{folg}}
\newcommand{\blem}{\begin{lem}}
\newcommand{\elem}{\end{lem}}
\newcommand{\bnote}{\begin{note}}
\newcommand{\enote}{\end{note}}
\newcommand{\bprf}{\begin{proof}}
\newcommand{\eprf}{\end{proof}}
\newcommand{\bd}{\begin{displaymath}}
\newcommand{\ed}{\end{displaymath}}
\newcommand{\be}{\begin{eqnarray*}}
\newcommand{\ee}{\end{eqnarray*}}
\newcommand{\eeqa}{\end{eqnarray}}
\newcommand{\beqa}{\begin{eqnarray}}
\newcommand{\bi}{\begin{itemize}}
\newcommand{\ei}{\end{itemize}}
\newcommand{\bnum}{\begin{enumerate}}
\newcommand{\enum}{\end{enumerate}}
\newcommand{\la}{\langle}
\newcommand{\ra}{\rangle}
\newcommand{\eps}{\epsilon}
\newcommand{\ve}{\varepsilon}
\newcommand{\vp}{\varphi}
\newcommand{\lra}{\longrightarrow}
\newcommand{\Lra}{\Leftrightarrow}
\newcommand{\Ra}{\Rightarrow}
\newcommand{\sub}{\subset}
\newcommand{\ems}{\emptyset}
\newcommand{\sms}{\setminus}
\newcommand{\ints}{\int\limits}
\newcommand{\sums}{\sum\limits}
\newcommand{\lims}{\lim\limits}
\newcommand{\bcup}{\bigcup\limits}
\newcommand{\bcap}{\bigcap\limits}
\newcommand{\beq}{\begin{equation}}
\newcommand{\eeq}{\end{equation}}
\newcommand{\einhalb}{\frac{1}{2}}
\newcommand{\rr}{\mathbb{R}}
\newcommand{\rn}{\mathbb{R}^n}
\newcommand{\ccc}{\mathbb{C}}
\newcommand{\cn}{\mathbb{C}^n}
\newcommand{\M}{{\cal M}}
\newcommand{\drehgleich}{\mbox{\begin{rotate}{90}$=$  \end{rotate}}}
\newcommand{\turngleich}{\mbox{\begin{turn}{90}$=$  \end{turn}}}
\newcommand{\turnsimeq}{\mbox{\begin{turn}{270}$\simeq$  \end{turn}}}
\newcommand{\vf}{\varphi}
\newcommand{\earr}{\end{array}\]}
\newcommand{\barr}{\[\begin{array}}
\newcommand{\bvec}{\left(\begin{array}{c}}
\newcommand{\evec}{\end{array}\right)}
\newcommand{\sumk}{\sum_{k=1}^n}
\newcommand{\sumi}{\sum_{i=1}^n}
\newcommand{\suml}{\sum_{l=1}^n}
\newcommand{\sumj}{\sum_{j=1}^n}
\newcommand{\sumij}{\sum_{i,j=1}^n}
\newcommand{\suminf}{\sum_{k=0}^\infty}
\newcommand{\inv}{\frac{1}}
\newcommand{\wzbw}{\hfill $\Box$\\[0.2cm]}
\newcommand{\lag}{\mathfrak{g}}
\newcommand{\lan}{\mathfrak{n}}
\newcommand{\lah}{\mathfrak{h}}
\newcommand{\laz}{\mathfrak{z}}
\newcommand{\+}{\oplus}
\newcommand{\x}{\times}
\newcommand{\lx}{\ltimes}
\newcommand{\rrn}{\mathbb{R}^n}
\newcommand{\laso}{\mathfrak{so}}
\newcommand{\lason}{\mathfrak{so}(n)}
\newcommand{\lagl}{\mathfrak{gl}}
\newcommand{\lasl}{\mathfrak{sl}}
\newcommand{\lasp}{\mathfrak{sp}}
\newcommand{\lasu}{\mathfrak{su}}
\newcommand{\w}{\omega}
\newcommand{\pmh}{{\cal P}(M,h)}
\newcommand{\s}{\sigma}
\newcommand{\deri}{\frac{\partial}}
\newcommand{\ddx}{\frac{\partial}{\partial x}}
\newcommand{\ddz}{\frac{\partial}{\partial z}}
\newcommand{\ddi}{\frac{\partial}{\partial y_i}}
\newcommand{\ddj}{\frac{\partial}{\partial y_j}}
\newcommand{\ddk}{\frac{\partial}{\partial y_k}}
\newcommand{\ddp}{\frac{\partial}{\partial p_i}}
\newcommand{\ddq}{\frac{\partial}{\partial q_i}}
\newcommand{\xz}{^{(x,z)}}
\newcommand{\mh}{(M,h)}
\newcommand{\wxz}{W_{(x,z)}}
\newcommand{\qmh}{{\cal Q}(M,h)}
\newcommand{\bbem}{\begin{bem}}
\newcommand{\ebem}{\end{bem}}
\newcommand{\bbez}{\begin{bez}}
\newcommand{\ebez}{\end{bez}}
\newcommand{\bbsp}{\begin{bsp}}
\newcommand{\ebsp}{\end{bsp}}
\newcommand{\pr}{pr_{\lason}}
\newcommand{\huts}{\hat{\s}}
\newcommand{\whut}{\w^{\huts}}
\newcommand{\bhg}{{\cal B}_H(\lag)}
\newcommand{\aaa}{\alpha}
\newcommand{\bb}{\beta}
\newcommand{\laa}{\mf{a}}
%newcommand{\alpha}{\alpha}
\newcommand{\lam}{\lambda}
\newcommand{\LL}{\Lambda}
\newcommand{\D}{\Delta}
\newcommand{\ß}{\beta}
\newcommand{\ä}{\alpha}
\newcommand{\W}{\Omega}
\newcommand{\esel}{\ensuremath{\mathfrak{sl}(2,\ccc)}}
\newcommand{\kg}{{\cal K}(\lag)}
\newcommand{\bg}{{\cal B}_h(\lag)}
\newcommand{\kk}{  \mathbb{K}}
\newcommand{\xy}{[x,y]}
\newcommand{\perdef}{$\stackrel{\text{\tiny def}}{\iff}$}
\newcommand{\eqdef}{\stackrel{\text{\tiny def}}{=}}
\newcommand{\lai}{\mf{i}}
\newcommand{\lar}{\mf{r}}
\newcommand{\Dim}{\mathsf{dim\ }}
\newcommand{\im}{\mathsf{im\ }}
\newcommand{\Ker}{\mathsf{ker\ }}
\newcommand{\trace}{\mathsf{trace\ }}
\newcommand{\grad}{\mathsf{grad}}
\newcommand{\lecturecount}{\begin{center}{\sf  (Lecture \Roman{lecturenr})}\end{center}\addtocontents{toc}{{\sf  (Lecture \Roman{lecturenr})}} \refstepcounter{lecturenr}
}
%% CHOL
\newcommand{\T}{{\cal T}}
\newcommand{\cur}{{\cal R}} 
\newcommand{\pd}{{\cal P}}  

\swapnumbers
\theoremstyle{definition}
\newtheorem{de}{Definition}[section]
\newtheorem{bem}[de]{Remark}
\newtheorem{bez}[de]{Notation}
\newtheorem{bsp}[de]{Example}
\theoremstyle{plain}
\newtheorem{lem}[de]{Lemma}
\newtheorem{satz}[de]{Proposition}
\newtheorem{folg}[de]{Corollary}
\newtheorem{theo}[de]{Theorem}

\bibliographystyle{alpha}

%\markright{\centerline{{\sf --- Preliminary version ---}}}

\title{Conformal holonomy of C-spaces, Ricci-flat,\\and Lorentzian manifolds}

\author{Thomas Leistner\\[.3cm]
{\small {\em Pure Mathematics - School of Mathematical Sciences}}\\
{\small {\em The University of Adelaide, SA 5005, Australia}}\\[.2cm]
{\small {\em phone: +61 (0)8 83033712,
fax: +61 (0)8 83033696}}\\
{\small {\em email:} {\tt tleistne@maths.adelaide.edu.au}}
}

\maketitle
 \begin{abstract}
The main result of this paper is that a Lorentzian manifold is locally conformally equivalent to a manifold with recurrent lightlike vector field and totally isotropic Ricci tensor if and only if its conformal tractor holonomy admits a 2-dimensional totally isotropic invariant subspace.
Furthermore, for semi-Riemannian manifolds of arbitrary signature we prove that the conformal holonomy algebra of a C-space is a Berger algebra. 
For Ricci-flat spaces we show how the conformal holonomy can be obtained by the holonomy of the ambient metric and get results   for Riemannian manifolds and plane waves.
\\[.3cm]
\end{abstract}

%\tableofcontents

\section{Introduction}
The question whether a semi-Riemannian manifold is conformally equivalent to an Einstein manifold was of particular interest in the last decade.
A tool in order to  solve this question is the so-called {\em tractor bundle} with its {\em tractor connection} over a conformal manifold, which was introduced by T. Y. Thomas \cite{thomas26}, \cite{thomas32}, further developed by T.N. Bailey, M.G. Eastwood and A.R. Gover \cite{bailey-eastwood-gover94} and \cite{eastwood95} and extensively treated in papers of A.R Gover and A. \v{C}ap \cite{cap/gover00},
\cite{cap/gover02}, \cite{cap/gover03} and \cite{cap02}.
Parallel sections with respect to this connection are in one-to-one correspondence with metrics in the conformal class which are Einstein metrics.
The quest for parallel sections suggests to study the holonomy group of the tractor connection. By doing this other structures beside parallel sections, such as invariant spaces or forms become of interest and the final aim might be to know all possible tractor holonomy groups and the corresponding structures on the manifold.

In the main result of this paper we deal with  the very special case where the conformal structure has Lorentzian signature and the holonomy of the tractor connection admits a 2-dimensional, totally isotropic invariant subspace, a case which cannot occur for a Riemannian conformal structure. We will prove the following theorem.
\btheo\label{theo}
Let $(M,[h])$ be a conformal manifold of Lorentzian signature, ${\cal T}$ its conformal tractor bundle and $D$ the conformal tractor connection. The holonomy group of $D$ admits a 2-dimensional totally isotropic invariant subspace if and only if $(M,h)$ is locally conformally equivalent to a Lorentzian manifold with recurrent lightlike vector field and totally isotropic Ricci tensor.
\etheo
Regarding the classification problem of conformal holonomies this treats  the `problematic case' in Lorentzian signature in the following sense. 
If the conformal holonomy group does not act irreducibly the following cases may occur:
\bnum
\item
There is a {\em one-dimensional invariant subspace}. In this case the manifold is conformally Einstein, with zero scalar curvature if the invariant subspace is isotropic or with non-zero scalar curvature otherwise.
\item
There is a  {\em non-degenerate invariant subspace of dimension greater than $1$}. In this case the manifold is conformally equivalent to a product of Einstein spaces with related scalar curvature, and the tractor holonomy of the product is the product of the tractor holonomies  of the factors. Since these are Einstein, their tractor holonomy equals to the metric holonomy of the ambient metric which is known by the Berger classification. For these facts see \cite{leitner04killing} or generalise the results of \cite{armstrong05} from the Riemannian to the non-degenerate Lorentzian case.
\item
There is  {\em 2-dimensional, totally isotropic invariant subspace}. This case, which cannot occur in Riemannian signature, is studied in the present article and treated by theorem \ref{theo}.
\enum
We recall the result of \cite{leitner04killing} because it contains a proposition which is related to the third case, i.e. to our result, but with stronger assumptions and stronger conclusions.
\btheo\cite[Proposition 3]{leitner04killing}\label{fstheo}
Let ${\cal T}$ be the tractor bundle of a simply-connected manifold with conformal structure of signature $(r,s)$. The holonomy group of the tractor connection fixes a decomposable $(p+1)$-form on a fibre of the ${\cal T}$, for $1\le p < r+s$ if and only if one of the following cases occurs:
\bnum
\item There is a product of Einstein metrics  $g_1$ and $g_2$ 
in the conformal class of dimension $p$  resp. $r+s-p$ the scalar curvatures of which are related in the following way
\[S_2 \ p(p-1)=- S_1 (n-p)(n-p-1).\]
If $S_1\not=0$, then the tractor holonomy fixes a non-degenerate invariant subspace. If $S_1=0$ it fixes a degenerate subspace of dimension $(p+1)$ with one lightlike direction.
\item
There is a metric in the conformal class with totally isotropic Ricci tensor and a parallel, totally isotropic $p$-form.
\enum In the second case the tractor holonomy fixes a totally isotropic subspace of dimension at least $2$.
 \etheo
 The important implication of this theorem is the `only if'-direction. Here the assumption of theorem \ref{fstheo} are stronger than in our theorem \ref{theo}.  Obviously, the assumption  that the tractor holonomy fixes a decomposable $(p+1)$-form implies the existence of an $(n-p+1)$-dimensional invariant subspace: if $\aaa$ is the fixed form, then $\{v\in  V| v \lrcorner \aaa=0\}$ is an invariant subspace.  But the converse is only true if the subspace is non-degenerate. To see this consider the $\rr^{n+4}$ with a inner product $\la.,.\ra$ of index $(2,n+2)$, given by the matrix
$\left (\begin{array}{ccccc}
0 &0 &0 & 0 & 1\\ 
0 &0 &0 & 1& 0\\ 
0 &0 & E_{n} & 0 & 0 \\ 
0 &1 &0 &0 & 0 \\ 
1 &0 &0 &0 & 0 \\ 
\end{array}\right)$ with respect to  the coordinates $(x_1,x_2,y_1, \ldots , y_n, z_1,z_2)$ on $\rr^{n+4}$, the $y_i$'s being spacelike and the remaining coordinates lightlike.
Consider as group $G$ the isotropy group of $L:=\mathsf{span} (x_1,x_2)$ in $SO(2,n+2)$ which has the Lie algebra
\[\lag:=\mf{iso}(L)
=\left\{\left. \left (\begin{array}{ccc}
X & \begin{array}{c}  u^t \\ v^t \end{array} &\begin{array}{rr}  a&0 \\ 0&-a \end{array}  \\
\begin{array}{cc} 0&0\end{array} &  A& \begin{array}{cc} -v&-u\end{array}\\
\begin{array}{cc}  0&0 \\ 
0&0 \end{array} & \begin{array}{c}  0 \\ 0 \end{array} & -X^t
\end{array}\right)\right|\, \begin{array}{l}
X\in\lagl(2),\\
 A\in\laso(n),\\ u,v\in\rr^n,\\c\in \rr\end{array}\right\}.\]
 By definition, $G$ has no other invariant subspace except $L$ and $L^\bot$, the latter spanned by $x_1,x_2, y_1, \ldots y_n$. 
If there is a decomposable $(n+2)$-form $\aaa$ such that $\lag\cdot \aaa=0$, then 
$L=\{v\in \rr^{n+4}| v\lrcorner \aaa=0\}$, which implies that $\aaa=a_1\la x_1.,\ra\wedge
a_2\la x_2,.\ra \wedge b_1\la y_1,.\ra \wedge \ldots \wedge b_n\la y_n,.\ra$.  But this form does not satisfy that $A\cdot \aaa=0$ for all $A\in \lag$, because 
\[\left(\left(\begin{array}{rcr}E_2&0&0\\0&0&0\\0&0&-E_2\end{array}\right)\cdot \aaa\right) (z_1,z_2, y_1, \ldots , y_n)=-a_1\cdot a_2\not=0.\]

In Lorentzian signature the second case of theorem \ref{fstheo} implies that the manifold is conformally equivalent to a Lorentzian manifold with parallel lightlike vector field, a so-called {\em Brinkmann wave} and with totally isotropic Ricci tensor, a conclusion which is obviously stronger than ours. But the above
  algebraic example shows that, in general, not only the conclusion but  also the assumption of theorem \ref{fstheo} is stronger.
  
 Since the equivalence `invariant form $\iff$ invariant subspace' fails in the case where the subspace is totally isotropic, we use a different approach than  the one of \cite{leitner04killing} which is based on this equivalence.
 
 \bigskip
 
The structure of the paper is as follows: in the introductory section we present the basic notions of conformal geometry such as density and tractor bundles, the tractor connection and its holonomy and recall their basic properties. For sake of brevity we ignore the relations to Cartan connections and parabolic geometries which can be found  in  \cite {kobayashi72}, \cite{cap/slovak/soucek97-2}, \cite{cap/slovak/soucek97-1} and \cite{cap/schichl00} and in various papers cited in  section 2.

Then we try, analogously as it was done  for metric holonomies, to derive algebraic constraints to the tractor holonomy based on the Bianchi-identity of the Weyl and the Schouten-Weyl tensor. We obtain the result that the tractor holonomy algebra is a {\em Berger algebra} if the conformal class contains the metric of a C-space. This result applies to the case where the conformal class contains a locally symmetric metric. 

Then we recall the fact that if the conformal class contains an Einstein metric $g$ the holonomy of the ambient metric of $g$  is equal to the tractor holonomy. We prove this fact in case of Ricci-flat manifolds: if the conformal class contains a Ricci-flat metric $g$, then its tractor holonomy equals to the holonomy of the ambient metric of $g$ which is contained in a semi-direct product of the holonomy of $g$ and $\rrn$. As a corollary we obtain a classification of `indecomposable Ricci-flat tractor holonomies' in the Riemannian case, which was independently obtained by \cite{armstrong05}.

In the next section we turn to Lorentzian manifolds with a recurrent lightlike vector field, describe their basic properties and  prove some results about their Ricci curvature which we will need in the proof of theorem \ref{theo}. As an interlude we  introduce `pr-waves' which are generalisations of pp-waves and calculate their metric holonomy. They provide an example of Lorentzian manifolds where the implication `recurrent lightlike vector field and totally isotropic Ricci tensor $\Rightarrow$ parallel lightlike vector field' is true.

Then we prove theorem \ref{theo}, using methods which are inspired by \cite{armstrong05}. Finally we calculate the conformal holonomy of plane waves and verify that they are conformally Ricci-flat.

\bigskip

{\bf Acknowledgements.} The author would like to thank Helga Baum and Michael Eastwood for their support and helpful discussions on the topic and Stuart Armstrong for making his work in progress \cite{armstrong05} available, which had great influence on the final form of some of the proofs.

\section{Conformal structures on manifolds}

\paragraph{Conformal changes of the metric.} 
Let $(M,g)$ be a semi-Riemannian manifold of dimension $n\ge 4$ with Levi-Civita connection $\nabla$. We will recall the definitions of the basic curvature quantities which are related to conformal properties of the manifold. 

$\cur(X,Y)=\left[\nabla_X,\nabla_Y\right] -\nabla_{[X,Y]}$ is the curvature tensor of $\nabla$,  $Ric=\trace_{(1,3)}\cur$ is the {\em Ricci tensor}, and $S=\trace Ric$ the {\em scalar curvature}. The trace adjusted Ricci tensor, or {\em Schouten tensor}, is defined as
\[P\ =\   \frac{1}{n-2} \left(   Ric- \frac{S}{2\left(n-1\right)}  g \right).\]
Its trace is equal to $\frac{S}{2(n-1)}$. The {\em Weyl tensor} is the traceless part of $\cur$ and given by
\[ W= {\cal R}-g \diamond P\]
where $A\diamond B$ is the  
the {\em Kulkarni-Numizu product} of two symmetric tensors $A$ and $B$:
\begin{eqnarray*}
(A\diamond B) (U,V,X,Y)& := &A(U,X) B(V,Y)+ A(V,Y) B(U,X) \\
&&
{}
 -A(U,Y) B(V,X)-A(V,X)B(U,Y)
\end{eqnarray*}
The {\em Schouten-Weyl tensor} is the skew-symmetrisation of $\nabla P$:
\[C(X,Y,Z)\ :=\ \left(\nabla_X P\right)(Y,Z) -  \left(\nabla_Y P\right)(X,Z).\]
It satisfies $(n-3)C(X,Y,Z)=(\mathsf{div}\  W)(X,Y,Z)=\trace_{(1,5)} (\nabla W)(X,Y,Z)$.

A manifold is {\em Einstein}, if $Ric = f \cdot g$ for a smooth function $f$. Tracing gives $f=\frac{S}{n}$ and the second Bianchi identity of $\cur$ implies that $S$ has to be constant.
The Schouten tensor of an Einstein manifold satisfies 
\beq\label{schouten}P\ =\ \frac{S}{2n(n-1)} \cdot g.\eeq
Furthermore, a manifold is called a {\em C-space} or {\em with harmonic Weyl tensor} if $C=0$. Of course, Einstein manifolds are C-spaces.

One is interested in finding the  conditions under which a metric is conformally equivalent to an Einstein metric.
Changing the metric $g$ conformally to a metric $\tilde{g}:=e^{2\vf}\cdot g$, where $\vf$ is a smooth function on $M$, one obtains the following transformation behavior 
\begin{eqnarray}
\widetilde{\nabla}_X Y&=& \nabla_X Y + d\vf(X) Y +d\vf(Y)X -g(X,Y)\mathsf{grad} \vf ,
\label{lctrafo}\\
\widetilde{\cur}&=& e^{2\vf}\left( \cur +g\diamond \Big(H_\vf-(d \vf)^2 +||\mathsf{grad}\vf||^2 \cdot g \Big)\right),\nonumber\\
\widetilde{P}&=&P - H_\vf +(d\vf)^2 -\einhalb ||\mathsf{grad}\vf||^2\cdot g,
\label{ptrafo}
\\
\widetilde{W}&=&e^{2\vf} W\nonumber,
\end{eqnarray}
where the quantities with \  $\widetilde{ }$\  \ are those of  the changed metric $\tilde{g}$, and $H_\vf$ is the symmetric {\em Hessian form} of $\vf$ defined by 
\[H_\vf(X,Y)\ =\ g(\nabla_X \mathsf{grad} \vf, Y)\ =\ (\nabla_X d\vf) (Y).\]
Equation (\ref{ptrafo}) shows that 
 $(M,g)$ is conformally Einstein if and only if $P - H_\vf +(d\vf)^2$ is a multiple of the metric, i.e. is pure trace.
By substituting $\vf=-\log \s$ for a non-vanishing function $\s$, i.e. $\tilde{g}=\s^{-2}\cdot g$, we get rid of the term $d\vf^2$ and this equation simplifies to
\beq H_\s + P \cdot \s \text{ is pure trace,}\label{einsteq}\eeq
because
$d\left(\log \s\right)=\frac{1}{\s} d\s$ and $H_{\log \s} = \frac{1}{\s}H_\s - \frac{1}{\s^2}d\s^2$.

\paragraph{Conformal structures and density bundles.} What is presented in this and the remaining paragraphs of this section mainly follows  \cite{bailey-eastwood-gover94}, \cite{eastwood95}, \cite{cap/gover00} and \cite{gover01}.

Let $(M,c)$ be a manifold of dimension $n$ with conformal structure $c$ of signature $(r,s)$. $c$ is an equivalence class of smooth semi-Riemannian metrics of signature $(r,s)$ which differ by a nowhere vanishing smooth function. 
The bundle of frames of a conformal manifold reduces to the bundle of frames which are orthonormal with respect to a metric from the conformal class $c$, denoted by 
${\cal CO}(M,c)=\{(e_1, \ldots , e_n)\in {\cal O}(M,g)|g\in c\}$.
The tangent bundle  can be associated to this bundle, 
$TM= {\cal CO}(M,c)\times_{CO(r,s)} \rrn$ where $CO(r,s)=\rr^+\times SO_0(r,s)\subset Gl(n)$ is the conformal group with its standard representation on $\rrn$. Its center $\rr^+$ acts on $CO(M,c)$ via $t \cdot (e_1, \ldots , e_n)= (t^{-1}e_1, \ldots ,t^{-1} e_n) $. The ray bundle  $\cal Q:=\bigcup_{p\in M}\{g_p|g\in c\}\subset \odot^2TM$ in the bundle of metrics of signature $(r,s)$ is a principle $\rr^+$-fibre bundle bundle. The group homomorphism
\barr{rcrcl}
\mathsf{det}^{2/n}&:& CO(r,s)&\rightarrow &\rr^+\\
&&A=a\cdot Id \times \tilde{A} &\mapsto & \mathsf{det}(A)^{2/n}=a^2.
\earr
and the bundle homomorphism
\barr{rcrcl}
f&:&CO(M)&\rightarrow & \cal Q \\
&&(e_1, \ldots , e_n)&\mapsto & g\ \text{ if } (e_1, \ldots , e_n)\in{\cal O}_p(M,g)
\earr
define a $\mathsf{det}^{2/n}$-reduction of $CO(M)$ to $\cal Q$, i.e. the following diagramm commutes:
\barr{ccccl}
 CO(r,s)\times CO(M)&\longrightarrow &CO(M)&&\\
&&&\searrow&\\
\makebox[0cm][r]{{\scriptsize  $\mathsf{det}^{n/2}\times f$}}\downarrow&\circlearrowright&
\makebox[0cm][r]{{\scriptsize $f$}}\downarrow\makebox[0cm][l]{ $\;\;\circlearrowright$}&&M.\\
&&&\nearrow&\\
\rr^+\times \cal Q&\longrightarrow &\cal Q&&
\earr

 A conformal structure can be described by  a section in a certain vector bundle which is related to {\em densitiy bundles}. Density bundles allow us to write quantities on a conformal manifold in an invariant manner.
Let $\delta^w$ be the  representation of weight $w$ on $\rr$, either of the conformal group, i.e. $\delta^w (A) :t\mapsto \mathsf{det}(A)^{\frac{w}{n}} \cdot t $, or of $\rr^+$, i.e. $ \delta^w(\alpha)t=\alpha^w\cdot t$. Then the density bundles are associated vector bundles to this representation,
\[\cal E[w]\ :=\   {\cal CO}(M)\times_{ \delta^w} \rr \ =\ {\cal Q}\times_{\delta^w}\rr ,\]
i.e. $[\vf^2g, \vf^w]=[g, 1]=[(e_1, \ldots , e_n), 1]=[(\vf^{-1}e_1, \dots, \vf^{-1}e_n),\vf^w]$.  Densities can be multiplied with each other, $\cal E[v]\otimes \cal E[w] =\cal E[v+w]$, and with the usual tensor bundles,
\[TM[w]:= \cal E[w]\otimes TM\ ,\  T^*M[w]:=\cal E[w]\otimes T^*M, \text{ etc. for tensor products of $TM$.}\]
If $\W\in \otimes^kT^*M[w]$ and $X_i\in TM[v_i]$ for $i=1,\ldots ,k$ one obtains $\W(X_1,\ldots , X_k)\in \cal E [w+v_1+ \ldots v_k]$.

Many important tensors can be considered as sections in these bundles. E.g. for a function $f\in C^\infty (M)$ its gradient w.r.t. to a metric in the conformal class becomes an invariant element in $TM[-2]$: $\mathsf{grad} f := [g,\grad^g(\vf)]=[\vf^2 g, \vf^{-2}\grad^g(f)]=[\vf^2g,\grad^{\vf^2g}(f)]$. The Weyltensor becomes an element in $\odot^2\Lambda^2T^*M[2]$. 
%Now, any metric in the conformal class $g\in c$ can be identified with a 
% density$\mu=[g,1]\in\cal E [-2]$. 
A conformal class $c$ can be seen as a section in  $\odot^2 T^*M[2]$:
\be
c&\mapsto& 
%[(e_1,\ldots , e_n), 1] \cdot g =
 [g,1] \otimes g \in\Gamma(\odot^2 T^*M[2]),
\ee
where $g$ is in the conformal class $c$.
%$(e_1,\ldots , e_n)\in {\cal O}(M,g)$ is a orthonormal frame with respect to a metric $g$ in the conformal class $g$. 
This map is well-defined as
$\vf^2 g\in c$ is mapped to 
%$[(\vf^{-1}e_1,\ldots , \vf^{-1} e_n), 1] \cdot 2\vf g=[(\vf^{-1}e_1,\ldots , \vf^{-1} e_n), 2\vf]\cdot g = [(e_1,\ldots , e_n), 1] \cdot g$ resp. to 
$[\vf^2g, 1]\otimes  \vf^2g=
 [\vf^2g, \vf^2]\otimes  g= [g,1]\otimes g$. On the other hand the conformal class $c=[g]$ provides a bundle homomorphism
\barr{rcccl}
[g]&:& TM[w]\times TM[v]&\rightarrow &\cal E[w+v+2]\\
 && (X,Y)= \left([g,1]\otimes X, [g,1]\otimes Y\right)&\mapsto & [g, g(X,Y)].\earr
Again,  this is well-defined as $[g]\left([\vf^2 g, \vf^w X], [\vf^2 g, \vf^v Y]\right)=[\vf^2 g, \vf^{w+v+2}g( X,Y)]$. For $v=w$ it is symmetric and for $v=w=-1$ it maps to $\cal E[0]=C^\infty (M)$, i.e. provides a proper metric on $TM[-1]$. As usual, $[g]$ identifies $TM[w]$ with $TM^*[w+2]$.
 
 Next, one defines the vector bundle
 \[\cal E\ := \  \cal E[1]\+TM[-1]\+\cal E[-1].\]
 It is equipped with metric $\la.,.\ra$ of signature $(r+1,s+1)$ by the formula
 \[\left\la (\s,X,\rho), (\xi,Y,\eta)\right\ra:= \s\eta +\rho \xi + [g](X,Y).\]
Any non-vanishing function $\vf$ on $M$ yields an bundle isomorphism $\Theta_\vf$ of $\cal E$  defined by the formula
 \[
 \Theta_\vf\left(\s,X,\rho\right):=\left(\s\ ,\  X+ \s\otimes \vf^{-1} \mathsf{grad}\ \vf\ ,\  
\rho -\vf^{-1}d\vf(X)-\einhalb  \vf^{-2} ||\mathsf{grad}\ \vf ||^ 2 \otimes \s \right),
\]
where $\grad\ \vf\in TM[-2]$ is the invariant gradient, $||\mathsf{grad}\ \vf ||^ 2=[g]\left( \mathsf{grad}\ \vf, \mathsf{grad}\ \vf\right)\in \cal E[-2]$, and the expression $d\vf(X)$ is meant to be $[g,X(\vf)]\in\cal E[-1]$ if $X$ is represented by $[g,X]\in TM[-1]$.
$\Theta_\vf$ is  an isometry with respect to $\la.,.\ra$, i.e. $\ \Theta_\vf|_p\in SO({\cal E}_p, \la.,.\ra_p)$.

\bigskip

Every metric $g$ in the conformal class defines a covariant derivative $D^g$ on the various densitiy bundles by
\[D^g_X[g,\vf]\ :=\  [g, X(\vf)].\]
Changing the metric in $c$ via $\tilde{g}=\vf^2 \cdot g$ it transforms as follows,
\beq
\label{denstrafo}
D^{\tilde{g}}_X\s\ =\ D^g_X\s +w\vf^{-1}X(\vf) \s, \text{ for }\s\in \Gamma(\cal E [w]),
\eeq
because for $\s=[g,f]$ it is
$
 D^{\tilde{g}}_X\s-D^g_X\s - =
[\vf^2g, X(\vf^w f)]-[g,X(f)]$.
This covariant derivative together with the Levi-Civita connection $\nabla^g$ of $g$ extends to covariant derivatives of the tensor bundles with densities $TM[w]$, $T^*M[w]$, etc., denoted by $D^g$ as well. They satisfy certain transformation formulae resulting from (\ref{denstrafo}) and (\ref{lctrafo}), e.g. for $Y\in TM[w]$ and $\w\in T^*M[w]$ it is
\barr{rcccl}
D^{\tilde{g}}_X Y&=&D^g_X Y+\vf^{-1}\left((w+1) X(\vf) \otimes Y +  Y(\vf) \otimes X-[g](X,Y)\otimes\mathsf{grad} ( \vf) \right)&\in&TM[w],\\
D^{\tilde{g}}_X \w&=&D^g_X \w+\vf^{-1}\left((w-1) X(\vf) \cdot \w -  \w(X)\otimes d\vf +\w\left(\mathsf{grad}( \vf)\right)\otimes [g](X,.) \right)&\in&T^*M[w].\earr
 
Finally, every metric defines a covariant derivative $D^g$ on $\cal E$  by the formula
\[D^g_X(\s,Y,\rho):=\left( D^g_X\s-[g](X,Y), D^g_X Y +\rho\otimes X+\s\otimes P^g(X)^\sharp, D^g_X\rho-P^g(X,Y)\right).\]
Here $P^g(X)^\sharp\in TM[-2]$ is the image of $P^g(X,.)$ under the identification of $T^*M$ with $TM[-2]$ via $[g]$. $P^g(X,Y)$ means $[g,P^g(X,Y)]\in \cal E[-1]$ for $Y$ being representated by $[g,Y]\in TM[-1]$.
$D^g$ is compatible with $\la.,.\ra$ and
 has the remarkable invariance property 
 \beq\label{dg}
 D^{\tilde{g}}_X\left(\Theta_\vf (\s,Y,\rho)\right)\ =\ \Theta_\vf \left(D^g_X(\s,Y,\rho)\right),
 \eeq
for $\tilde{g}=\vf^2\cdot g$.

\paragraph{Tractor bundle and tractor connection.}
The  transformation formulae above imply that the equation
\beq\label{confeinst}\text{trace-free part of }\left(D^g D^g +P^g\right)\s =0\eeq as a differential equation on $\Gamma({\cal E}[1])$ is conformally invariant. By equation (\ref{einsteq}) there is a solution  $\s=[g,\vf]$ of (\ref{confeinst}) if and only if the metric $\vf^{-2} g$ is an Einstein metric.
But considering the equation as an equation on the $2$-jet bundle of ${\cal  E}[1]$, denoted by $J^2  {\cal  E}[1]$, it has always solutions. The tractor bundle ${\cal T}$ is defined as its solution space in  $J^2  {\cal  E}[1]$.

\bigskip

Every metric $g$ in the conformal class together with its Levi-Civita connection $\nabla$  gives a more manageable realisation of the tractor bundle $\cal T$:  $g\in c$ defines an isomorphism
$\Psi^g : \cal T\ \rightarrow \cal E,$
such that the following diagram commutes
\barr{rrl}
&{ }_{\Psi^g}&\cal E\\[-.1cm]
&\;\nearrow&\\
{\cal T} &\circlearrowright&\Big.\Big\downarrow\ \text{{\scriptsize $\Theta_\vf$}}\\
&\;\searrow&\\[-.1cm]
&{}^{\Psi^{\tilde{g}}}\ &\cal E
\end{array},\ \text{
where $\tilde{g}=\vf^2\cdot g$.}\]
The invariance properties of $\la.,.\ra$ and $D^g$ ensure that both can be transferred to the tractor bundle, giving a metric and a compatible covariant derivative $D$ on $\cal T$.

We will denote by $\cal F$ the curvature of $D$ and by $Hol_p(M,c):=Hol_p(\cal T, D)$ the holonomy group of $D$, which is called {\em conformal holonomy} of $(M,c)$. Since $D$ is compatible with $\la.,.\ra$ the holonomy group is contained in $SO(r+1,s+1)$.
In the following calculations we will always fix a metric $g$ in the conformal class which identifies on one hand $\cal T$ with $\cal E$ and on the other hand $\cal E[\pm 1]$ with $C^\infty(M)$ and $TM[-1]$ with $TM$.  $D$  can be written as
\[D_X(\s,Y,\rho)\ =\ \left(X(\s)-g(X,Y)\ ,\; \nabla_X Y+\rho \cdot X+\s P(X)^\sharp \ ,\;X(\rho)-P(X,Y)\right),\]
where $g(P(X)^\sharp,.)=P(X,.)$.
We can express the curvature of $D$ as follows
\[
\cal F (X,Y)(\s,Z,\rho)\ =\ \left( 0\ ,\; \s C(X,Y)^\sharp + W(X,Y)Z\ ,\;-C(X,Y,Z)\right),\]
or in matrices
\[\cal F (X,Y)\ =\ \left(\begin{array}{ccc}
0&0&0\\[.2cm]
C(X,Y)^\sharp&W(X,Y)&0\\[.2cm]
0&-C(X,Y,Z)&0\end{array}\right).\]
Now we turn to important properties of the connection $D$ and how it is related to conformally Einstein metrics.
First we note the following fact for recurrent sections.
\blem\label{clem}
If $(\s,Y,\rho)\in \Gamma(\cal E )$ is a non-trivial {\em recurrent} section of $D$, i.e.
\beq\label{recur}D_X(\s,Y,\rho)\ = \ \theta(X)\cdot (\s,Y,\rho),\eeq
 then there is no open set on which 
$\s$ is zero, and sections of the form $(0,0,\rho)$ cannot be recurrent. Any recurrent section 
can locally be rescaled such that the rescaled section is parallel.  \elem
 \bprf
Suppose 
that $\s\equiv 0$ on an open subset. Then (\ref{recur}) gives that $g(X,Y)=0$ for all $X\in TM$ and thus $Y=0$ on this subset which implies that $\rho\cdot X=0$, i.e. $\rho=0$.  On  the other hand, $D_X(0,0,\rho)=(0,\rho X, X(\rho))$ cannot be a multiple of $(0,0,\rho)$.

For sections with non-zero length  the second point is obvious because dividing by the length always gives a parallel section. But for an isotropic recurrent section the statement is a special property of the tractor connection. Because of the first statement of the lemma we may assume that $\s\equiv 1 $ on an open subset. Then (\ref{recur}) implies that 
%\beq\label{recur1}
%\left(-g(X,Y), \nabla_X Y +\rho (X) + P(X)^\ast, X(\rho) -P(X,Y)\right)\ =\ \theta(X) (1,Y,\rho)\eeq
%i.e. 
\beqa\nonumber
\theta &=& -g(Y,.)\ \text{ and }\\
\label{recur2}
g(\nabla_U Y,V)&=& -\left(\rho g(U,V) + P(U,V)+ g(Y,U)g(Y,V)\right).
\eeqa

The rescaled section $f\cdot (1,Y,\rho)$ is parallel if $df=-\theta =g(Y,.)$, i.e.  we have to verify that the form $g(Y.,)$ is closed. But  (\ref{recur2}) ensures that
$\ d(g(Y.,))(U,V)\ =\ g(\nabla_U Y, V) - g(\nabla_VY,U)\ =\ 0$.
\eprf
The important property of the connection $D$ is that its parallel sections provide functions on a dense subset of $M$  which satisfy the conformally Einstein equation $(\ref{einsteq})$, i.e. correspond to Einstein metrics which are conformally equivalent to $g$: a parallel section $(\s,Y,\rho)$ satisfies 
\beqa
Y&=&\mathsf{grad}\ \s,
\nonumber
%\label{pars}
\\
\text{thus }\ \ 0&=&\ H_\s+\s P+\rho g,\label{pary}\\
\text{and by tracing }\ \  \rho&=&-\frac{1}{n} \left(\trace H_\s + \frac{S}{2(n-1)} \cdot \s\right).\nonumber
%\label{parrho}
\eeqa
By (\ref{pary}), $H_\s+\s P$  is pure trace, i.e. $\s^{-2}\cdot g$ is an Einstein metric on 
$M\setminus\{\s=0\}$. 
On the other hand, if we realise the tractor bundle and connetion by an Einstein metric $g$, then the section $(1,0, -\frac{S}{2n(n-1)})$, the length of which is $  -\frac{S}{2n(n-1)}$ is parallel by equation (\ref{schouten}).

Concluding this introductory remarks we recall that the definition of the tractor connection was independent of the chosen metric from the conformal class in order to obtain the following correspondences:
\be
&\left\{\begin{array}{c}
\text{Einstein metrics in the conformal class $c$ of a dense subset of $M\}$}\\ \text{with $S>0$ / $S<0$ / $S=0$}\end{array}\right\}&\\[.2cm]
&\updownarrow&\\[.2cm]
&\left\{
\begin{array}{c}
\text{parallel sections of the tractor connection}\\ \text{which are time- / space- / lightlike}\end{array}\right\}&\\[.2cm]
&\updownarrow&\\[.2cm]
&\left\{\begin{array}{c}
\text{one-dimensional $Hol_p(M,c)$-invariant subspaces}\\
\text{which are time- / space- / lightlike.}\end{array}\right\}.&\ee
Finally we would like to add some references of recent papers which deal with this correspondence and with tractor holonomy \cite{leitner04groups}, \cite{gover/nurowski04}, \cite{nurowski04} and \cite{gover04}.

\section{Algebraic constraints on the conformal holonomy of C-spaces}

In this section we will make a little algebraic step towards the answer of the question when a conformal holonomy is the holonomy of a  Levi-Civita connection.
Due to the Ambrose-Singer holonomy theorem \cite{as} and the Bianchi-identity of the curvature, every  holonomy algebra of a torsion-free affine connection  satisfies the algebraic criterion to be a {\em Berger algebra}. This notion is defined as follows: for  a linear Lie algebra
 $\lag\subset \lagl(E)$ one has the space of algebraic curvature endomorphisms
\be
%\label{kg}
{\cal K}({\mathfrak g})& := &  \{ R\in \Lambda^2 E^* \otimes {\mathfrak g} \ |\ R(x,y)z +
R(y,z)x + R(z,x)y=0\},\\
\text{ and\  }\ \ 
\underline{{\mathfrak g}}& :=& span \{ R(x,y)\ |\  x,y\in E, R\in {\cal K}({\mathfrak g})\}.
\ee
$\lag$ is called {\em Berger algebra} if $\lag=\underline{\lag}$.
Holonomy groups of linear torsion-free connections are Berger algebras and the irreducible amongst them were recently classified by \cite{schwachhoefer1} and \cite{schwachhoefer2} extending the classification of \cite{berger55}.
For the proof of our result we need a lemma about the tractor curvature.
\blem
The curvature of the conformal tractor connection satisfies the following identity
\be
{\cal F}(X_1,X_2)(s_3,X_3,r_3)+
{\cal F}(X_2,X_3)(s_1,X_1,r_1)
+{\cal F}(X_3,X_1)(s_2,X_2,r_2)&
=&\\
\left( 0\ ,\  s_1\cdot C(X_2,X_3)^\sharp + s_2\cdot C(X_3,X_1)^\sharp +s_3\cdot C(X_1,X_2)^\sharp\ ,\ 0\right).\ee
\elem
\bprf
This equation follows immediately from the Bianchi identity for the Weyl tensor and the Schouten-Weyl tensor.
\eprf  
%This lemma implies an algebraic constraint for for the conformal holonomy in the case where the conformal class contains a C-space metric.
%\blem
%Let $(M,c)$ be a semi-Riemannian conformal manifold such that there is a metric $g$ in the conformal class $c$ such that $g$ si a metric of a C-space, i.e. the Schouten-Weyl tensor $C$ vanishes. Then the conformal curvature in a point $p$ satisfies
%\[

\btheo\label{ctheo}
Let $(M,c)$ be a  conformal manifold of arbitrary signature. If the conformal class $c$ contains the metric of a a C-space, then its conformal holonomy algebra is a Berger algebra.
\etheo

\bprf
Suppose that $g$ is the metric in the conformal class $c$ which has the property that its Schouten-Weyl tensor $C:=C^g$ vanishes. We consider the splitting of the tractor bundle and the formula for the tractor connection with respect to the metric $g$. Since $C=0$ we obtain for the tractor curvature
\[{\cal F}(X,Y)(s,Z,r)=(0,W(X,Y)Z,0)\]
for $X,Y,Z\in TM$, $s,r\in \rr$, i.e. $(s,Z,r)\in \T$ and $W $ the Weyl tensor with respect to the metric $g$.

Lets denote by $\pd_\gamma $ the parallel displacement with respect to the tractor connection $D$ along a curve $\gamma$. Corresponding to the decomposition of the tractor bundle according to the metric $g$ we may split this parallel displacement into components:
\be
\pd_\gamma&=& (\pd^-_\gamma, \pd^0_\gamma , \pd^+_\gamma)\in End (\T_{\gamma(0)}, \T_{\gamma(1)} ).
%,\text{ with}\\
%\pd^-_\gamma&:&{\cal E}_p[1]\+T_pM[1]\+{\cal E}_p[-1]\rightarrow {\cal E}_p[1]\\
%\pd^0_\gamma&:&{\cal E}_p[1]\+T_pM[1]\+{\cal E}_p[-1]\rightarrow TM_p[1]\\
%\pd^-_\gamma&:&{\cal E}_p[1]\+T_pM[1]\+{\cal E}_p[-1]\rightarrow {\cal E}_p[-1]
\ee
Since $\pd^0_\gamma: \T_{\gamma(0)}\rightarrow T_{\gamma(1)} M$ is surjective, 
 by the Ambrose-Singer holonomy theorem \cite{as} the holonomy algebra of the conformal connnection $D$ is given by
\be
\lefteqn{\mf{hol}_p(M,[h])=}\\
%=\mf{hol}_p(\T,D)
&=&\left\{\left. \left(\pd_{\gamma} \right)^{-1}\circ {\cal F}(X,Y)\circ  \pd_{\gamma} \right| \gamma \text{ a curve starting at $p$ and }X,Y\in T_{\gamma(1)}M\right\}\\
&=&\left\{\left. \left(\pd_{\gamma}^\T \right)^{-1}\circ {\cal F}\left(\pd^0_\gamma(s,X,r),\pd^0_\gamma (u,Y,v)\right)\circ  \pd_{\gamma} \right| 
(s,X,r), (u,Y,v)
\in \T_{\gamma(1)}\right\}.
\ee
Since $g$ is the metric of a C-space, i.e. $C^g=0$,  we obtain by the previous lemma
\[
\left(\pd_{\gamma}^\T \right)^{-1}\circ {\cal F}\left(\pd^0_\gamma(\ .\ ),\pd^0_\gamma (\ .\ )\right)\circ  \pd_{\gamma} \in {\cal K}\left(\mf{hol}_p(M,c)\right).\]
Hence, $\mf{hol}_p(M,c)\subset \laso( \T_p, \la.,\ra_p)$ is a Berger algebra.
\eprf

This result immediately  gives a consequence for locally symmetric spaces.
\bfolg
Let $(M,g)$ be a semi-Riemannian manifold which is locally conformally equivalent to a locally symmetric space. Then its conformal holonomy algebra is a Berger algebra.
\efolg

\bprf
For symmetric spaces the derivatives of curvature tensors vanish, i.e. $C=0$.
\eprf

%Furthermore it enables us to conclude in the Riemannian case that the conformal holonomy is the metric holonomy of a Lorentzian manifold, because Berger algebras of signatrue $(1,n+1)$ wher
% 
% Folgerung fuer symmetrische Raeume

\section{Ambient metrics for conformally Einstein spaces}

The idea of using  an ambient metric in order to describe conformal structures goes back to C. Fefferman and C. Robin Graham \cite{fefferman/graham85}, for recent results see \cite{fefferman/graham02}, \cite{cap02} \cite{fefferman/hirachi03}, 
\cite{cap/gover03}, \cite{graham/hirach04} and \cite{fefferman/hirachi03}.
It is known that if a semi-Riemannian manifold $(M,g)$ of signature $(p,q)$ is conformally Einstein, then there exists a semi-Riemannian manifold $(\overline{M},\bar{g})$ of signature $(p+1,q+1)$ which admits a parallel vector field the length of which depends on the sign of the scalar curvature of the Einstein metric in the conformal class of $g$.  The holonomy of $(\overline{M}, \bar{g})$ then is equal to the conformal holonomy of $(M,[g])$. We will recall this result and prove it in the case where we could find no proof in the literature.

\bs\cite{leitner04killing}
Let $(M,g)$ be a semi-Riemannian manifold of signature $(p,q)$ and of dimension $n:=p+q$. If $(M,g)$ is an Einstein space of non-zero scalar curvature $S$, then  the manifold $\overline{M}:=\rr\times M\times \rr^+$ with semi-Riemannian metric 
\[\bar{g}:=\frac{n(n-1)}{S}\left(dt^2-ds^2\right)+t^2 \cdot g\] has the property that  the holonomy of the Levi-Civita connection of $(\overline{M},\bar{g})$ and the  conformal holonomy of $(M,[g])$ coincide:
$Hol_{(1,p,1)}(\overline{M},\bar{g})\ =\ Hol_p(M,[g])$.
\es
The proof  in \cite{leitner04killing} relies on the decomposition of the tractor bundle with respect to the Einstein metric $g$ and on the following identification of the tangent bundle of $(\overline{M},\bar{g})$ with the tractor bundle of $(M,[g])$:
\barr{rcccccl}
\T\supset TM&\ni& (0,X,0)& \mapsto & X&\in& TM\subset T
\overline{M}|_{\{1\}\times M\times\{1\}}\\
\T&\ni & (1,0,\frac{S}{2n(n-1)})&\mapsto& \frac{S}{n(n-1)} \frac{\partial}{\partial t}&\in& T\overline{M}|_{\{1\}\times M\times\{1\}}\\
\T&\ni & (1,0,-\frac{S}{2n(n-1)})&\mapsto& \frac{S}{n(n-1)} \frac{\partial}{\partial s}&\in &T\overline{M}|_{\{1\}\times M\times\{1\}}.
\earr
The parallel vector field of the ambient metric $\bar{g}$ is equal  to $\frac{\partial}{\partial s}$ which is spacelike if $S<0$ and timelike if $S>0$.
Since the manifold $(\overline{M}, \overline{g})$ has a parallel vector field of non-zero length, its holonomy is equal to the holonomy of the cone 
\[\left(\ \hat{M}:= \rr^+\times M\ ,\  \hat{g}= \frac{n(n-1)}{S} dt^2 +t^2 g\ \right).\]
For Riemannian Einstein metrics $g$ with $S>0$ this cone is a Riemannian manifold. Assuming that $(M,g)$ is complete, the cone is either flat --- i.e. $(M,g)$ is locally isometric to the sphere and thus conformally flat --- or irreducible \cite{gallot79}. By the O'Neill formulas the cone over an Einstein manifold is Ricci flat \cite{oneill66}. This restricts the holonomy of the cone further and we obtain by the Berger list that it equals to $SO(n+1),\ SU(m)$ if $2m=n+1$, $ Sp(m)$ if $4m=n+1$, $G_2$ if $n=6$ or $Spin(7)$ if $n=7$.
If $S<0$ this cone is a Lorentzian manifold whose holonomy can be obtained by the Berger list and the classification of indecomposable, non-irreducible Lorentzian holonomy groups in
\cite{leistner02}, \cite{leistner03} and \cite{leistner03b}.

Now we turn to the case where the manifold is conformally Ricci-flat. Here a similar result holds. First we prove a general theorem about metric holonomy.
\btheo
Let $(M,g)$ be a semi-Riemannian manifold of signature $(r,s)$ and of dimension $n:=r+s$. Then the holonomy
group  of the pseudo-Riemannian manifold 
\[\left(\ 
\overline{M}:=\rr\times M\times \rr^+\ ,\ 
\bar{g}:=2dxdz+z^2\cdot g\ \right)\] 
 in the point $(1,p,1)$ is generated only by curves which
 run in $\{1\}\times M \times \{1\}$. If $k$ is the number of linear independent parallel vector fields on $(M,g)$, it satisfies
 \[ Hol_{(1,p,1)}(\overline{M},\overline{g})\ \subset \  Hol_p(M,g)\ltimes \rr^{n-k}\
 =\ \left\{\left.\left( \begin{array}{rrr}1 & v^t &-\einhalb v^t  v\\
 0&A&-Av\\
 0&0&1\end{array}\right)\right| \begin{array}{l} A\in Hol_p(M,h),\\ v\in \rr^{n-k}\end{array}\right\},
 \]
and 
$pr_{SO(n)} Hol_{(1,p,1)}(\overline{M},\overline{g})= Hol_p(M,g)$.
 \etheo
\bprf
Since $X=\ddx$ is a parallel vector field of $\left(\overline{M},\overline{g}\right)$, we notice that
\[Hol_{(1,p,1)}(\overline{M},\overline{g})\ \subset\  SO(r,s)\ltimes \rr^n\ \subset \ SO(T_{(1,p,1)}\overline{M},\overline{g}_{(1,p,1)}).\]
  Fixing coordinates $(y_1, \ldots , y_n)$ on $M$ we obtain as
the non-vanishing Christoffel symbols of the new Levi-Civita connection $\overline{\nabla}$ only the following
 \be
 \overline{\Gamma}_{ij}^0(x,p,z)&=&-z g_{ij}(p), \\
 \overline{\Gamma}_{ij}^k(x,p,z)&=&  \Gamma_{ij}^k(p)\text{ and}\\
 \overline{\Gamma}_{i\ n+1}^j(x,p,z)&=&\frac{1}{z}\delta_{ij},\ee
 where the indices $0$ and $n+1$ refer to the $x$- and $z$-coordinate.
 In other words, the non-vanishing covariant derivatives are the following
 \beqa\label{formeli}
 \overline{\nabla}_{Y_i}{Y_j}&=& \nabla_{Y_i}{Y_j} - z g_{ij} X\text{ and}\\
\label{formelz} \overline{\nabla}_{Y_i}{Z}= \overline{\nabla}_Z{Y_i}&=& \frac{1}{z} Y_i,
 \eeqa
 where $X=\ddx$, $Z=\ddz$ and $Y_i=\ddi$.
  First we take  a curve $\gamma(t)=\left(\gamma_1(t), \ldots , \gamma_n(t)\right) $ running in $\{x_0\}\times M \times \{z_0\}$ and \[U(t)=a(t) X(\gamma(t))+ Y(t)+c(t)Z(\gamma(t))\] the parallel displacement with respect to $\overline{\nabla}$ along $\gamma$ where $Y(t)=\sumk b_k(t) Y_k(\gamma(t))$ is the component in $\{x_0\}\times M \times \{z_0\}$. $U(t)$ satisfies
 \be
 0
 &=&
 \overline{\nabla}_{\dot{\gamma}(t)} U(t)
 \\
 &=&
\big( \dot{a}(t)- z_0 g(\dot{\gamma}(t), Y(t))\big) \cdot X(\gamma(t))
+\nabla_{\dot{\gamma}(t)}Y(t)
 + \frac{c(t)}{z_0}\dot{\gamma}(t)+\dot{c}(t) Z(\gamma(t)).
\ee
Hence, $\dot{c}(t)\equiv 0$.
If we assume that $U(0)\in T_{(x_0,p,z_0)}M$, i.e. $a(0)=c(0)=0$, then $c\equiv 0$ and thus
\beq\label{par1}
0\ =\ \Big( \dot{a}(t)- z_0 g(\dot{\gamma}(t), Y(t))\Big) X(\gamma(t)) +\nabla_{\dot{\gamma}(t)}Y(t).
\eeq
This implies that $Y(t)$ has to be the parallel displacement of $U(0)=Y(0)$ with respect to $\nabla$, i.e.
\beq\label{par2}
0=\sumij\dot{b}_k(t) +b_i\ \dot{\gamma}_j\Gamma^k_{ij} (\gamma(t))\text{ for all }k,\eeq
and that $a$ satisfies
\beq\label{par3}
\dot{a}(t)=z_0g\big(\dot{\gamma}(t), Y(t)\big) = z_0\sumk b_k\  g\big(\dot{\gamma}(t), Y_k(\gamma(t)))\big).\eeq
This implies that 
$Hol_p(M,g)\subset pr_{SO(n)}\left(Hol_{(1,p,1)}(\overline{M},\overline{g})\right)  $.

Now we consider a general curve $\overline{\gamma}(t)=\left( x(t), \gamma(t), z(t)\right)$ with $z(0)=z_0$ and the following vector field along $\overline{\gamma}(t)$:
\be
W(t)&=&\frac{z_0}{z(t)}Y(\overline{\gamma}(t)) +a(t) X(\overline{\gamma}(t))\\
&=&
\frac{z_0}{z(t)}
\sumk b_k(t) Y_k(\overline{\gamma}(t))+a(t) X(\overline{\gamma}(t)),
\ee
where $Y(t)$ is the parallel displacement of $Y(0)$ along $\gamma$ with respect to $\nabla$, i.e. $b_k$ satisfies (\ref{par2}) and $a$ satisfies (\ref{par3}) with respect to $\gamma(t)$ and $Y(t)$.
Then $W(t)$ is parallel along $\overline{\gamma}(t)$. In detail:
\be
\overline{\nabla}_{\dot{\overline{\gamma}}(t)}W(t)
&=&
-z_0\frac{\dot{z}(t)}{z^2(t)}\left( \sumk b_k(t) Y_k(t)\right)
+
\frac{z_0}{z(t)}\left( \sumk \dot{b}_k(t) Y_k(t)\right) \makebox[3cm][l]{}
\\
&&+
\frac{z_0}{z(t)} \sumk {b}_k(t) 
\underbrace{\overline{\nabla}_{\dot{\overline{\gamma}}(t)}Y_k(t)}_{
\begin{minipage}[b]{3cm}{\small \vspace{-.3cm}
\barr{cl}
=&\sumij \dot{\gamma}_j(t)\overline{\Gamma}^i_{jk}(\overline{\gamma}(t))Y_i(t)
\\[.2cm]
&
- z(t)\left(\sumi \dot{\gamma}_i(t)g_{ij}(\gamma(t)) \right)X(\overline{\gamma}(t))
\\[.2cm]
&
+
\frac{\dot{z}(t)}{z(t)}Y_k(t)
\earr
}
\end{minipage}
}
\\
&&
+\ 
\dot{a}(t)X(\overline{\gamma}(t)) + a(t)\underbrace{ \overline{\nabla}_{\dot{\Gamma}(t)}X(\overline{\gamma}(t))}_{=0}\\
&=&
\frac{z_0}{z(t)}\left( \sumk \dot{b}_k(t) +
b_i(t) \ \dot{\gamma}_j(t)\underbrace{\overline{\Gamma}^k_{ij}(\overline{\gamma}(t))}_{=\Gamma^k_{ij}(\gamma(t)}
\right)
Y_k(t)\\&&
+\left( \dot{a}(t)- z_0 \sum_{k,l=1}^n b_k(t)\dot{\gamma}_l(t) g_{kl}(\gamma(t))\right) X(\overline{\gamma}(t))
\\&=&0
\ee
by (\ref{par2}) und (\ref{par3}).
Since $Y(t)$ was the parallel displacement along $\gamma$ with respect to $\nabla$ the shape of $W(t)$ shows that $Hol_{(1,p,1)}(\overline{M},\overline{g})$ is generated only by curves which run in $\{1\}\times M\times \{1\}$ and is therefore contained in $  Hol_p (M,g)\ltimes \rrn$.

Finally, if $Y$ is a parallel vector field on $(M,g)$ the vector field $\frac{1}{z}Y$ is parallel on $(\overline{M}, \overline{g})$ and has values in $TM\subset T\overline{M}$. If there are $k$ parallel vector fields on $(M,g)$ this implies that $Hol_{(1,p,1)}(\overline{M},\overline{g}) \subset  Hol_p (M,g)\ltimes \rr^{n-k}$.
 \eprf
Before we draw the consequences for the tractor holonomy we will give a more explicit description of the ambient holonomy algebra, which we will need later on. First one calculates that the curvature of $\overline{g}$ reduces to the one of $g$:
\beqa\label{formelcur}
\overline{\cur}(U,V)&=&
\left(\begin{array}{ccc}
0&0&0\\[.2cm]
0&\cur(U,V)&0\\[.2cm]
0&0&0\end{array}\right) \text{ if } U,V\in TM,\text{ and $0$ otherwise.}
\eeqa
The holonomy algebra $\mf{hol}_{(1,p,1)}(\ol{M},\ol{g})$ is generated by expressions of the form
$\overline{\cal P}^{-1}_{\gamma(1)} \circ \ol{\cal R}_{\gamma(1)}(U,V)\circ \overline{\cal P}_{\gamma(1)}$ where $\overline{\cal P}_{\gamma(1)}$ denotes the parallel displacement w.r.t. $\ol{g}$  along curves $\gamma$ from $(1,p,1)$ to $\gamma(1)$, and $U,V\in T_{\gamma(1)}\ol{M}$. By the formulas above these terms can be calculated as follows with respect to the frame field $(X, Y_1, \ldots , Y_k, Z)$, for a curve $\ol{\gamma}(t)=(x(t),\gamma(t),z(t))$, $U,V\in T_{\ol{\gamma}(1)} M$,  and $\cal P_{\gamma(1)}$ the parallel displacement w.r.t. $g$:
\begin{eqnarray}
 \nonumber
\overline{\cal P}^{-1}_{\ol{\gamma}(1)} \circ \ol{\cal R}_{\ol{\gamma}(1)}(U,V)\circ \overline{\cal P}_{\ol{\gamma}(1)}(Y_i)
&=&
\frac{1}{z}\overline{\cal P}^{-1}_{\ol{\gamma}(1)}\left( {\cal R}_{\gamma(1)}(U,V) {\cal P}_{\gamma(1)}(Y_i)\right)\\
&=&a_i(1) X+  {\cal P}^{-1}_{\gamma(1)}\circ {\cal R}_{\gamma(1)}(U,V) \circ {\cal P}_{\gamma(1)}(Y_i)\label{pareqn}
\end{eqnarray}
where $a_i$ is determined by the following equation obtained by (\ref{par3}) in which $\gamma^{-1}(t)=\gamma(1-t)$
\begin{eqnarray}
 \nonumber
\dot{a}_i(t)&=& 
g\left(\dot{\gamma}^{-1}(t), \cal P_{\gamma^{-1}|_{[0,t]}}\cal R_{\gamma(1)}(U,V)\cal P_{\gamma(1)}(Y_i)\right)
\\
&=&-g\left(\dot{\gamma}(1-t), \cal P_{\gamma|_{[0,1-t]}}\left( \cal P_{\gamma(1)}^{-1}\circ\cal R_{\gamma(1)}(U,V)\circ \cal P_{\gamma(1)}(Y_i)\right)\right)\label{aeqn}
\end{eqnarray}
 
 \bfolg\label{richol}
 Let $(M,[g])$ be a conformal manifold of arbitrary signature where $g$ is a Ricci-flat metric, $(\overline{M},\overline{g})$ constructed as above. Then the conformal holonomy $ Hol_p(M,[g])$ is a pseudo-Riemannian holonomy and satisfies
\[ Hol_p(M,[g])\ =\ Hol_{(1,p,1)}(\overline{M},\overline{g})\ \subset\ 
Hol_p(M,g)\ltimes \rr^{n-k}.\]
\efolg
\bprf
By the identification $\Psi$ of the decomposed tractor bundle of $(M,g)$ and $T\overline{M}$ given by
\barr{rcccccl}
\T\supset TM&\ni& (0,Y,0)& \stackrel{\Psi}{\mapsto} & Y&\in& TM\subset T\overline{M}|_{\{1\}\times M\times\{1\}}\\
\T&\ni & (1,0,0)& \stackrel{\Psi}{\mapsto}&  \frac{\partial}{\partial x}&\in& T\overline{M}|_{\{1\}\times M\times\{1\}}\\
\T&\ni & (0,0,1)&\stackrel{\Psi}{\mapsto}&  \frac{\partial}{\partial z}&\in &T\overline{M}|_{\{1\}\times M\times\{1\}}
\earr
and the previous theorem we obtain the statement.
\eprf
In Riemannian signature we obtain the following result which was proven by \cite{armstrong05} without using the ambient construction. 
\bfolg
A Riemannian manifold $(M,g)$ is conformally Ricci-flat with indecomposable ambient metric if and only if its conformal holonomy is equal to the Lorentzian holonomy 
\[Hol_p(M,g)\ltimes \rr^{n}\]
of  a Brinkmann wave where $Hol_p(M,g)$ is a product of
 $SO(n),\ 
SU(m),\
Sp(m),\
Spin(7)$ or 
$G_2$.
\efolg
\bprf
Since the ambient metric is supposed to be indecomposable (in the sense of the next section) the projection of the tractor holonomy onto $\rrn$ is the whole $\rrn$. Hence $(M,g)$ has no parallel vector fields, i.e. $Hol(M,g)$ decomposes into irreducible componenents. By the description of indecomposable, non-irreducible subalgebras of $\laso(1,n+1)$ in \cite{bb-ike93} this implies that the holonomy of the ambient metric is not of coupled type, which implies the statement.\eprf

\section{Lorentzian manifolds with recurrent lightlike vector field}
A vector field $X$ is called recurrent if $\nabla X = \Theta \otimes X$ where 
$\Theta $ is a one-form on $M$. If the length of a recurrent vector field is non-zero, or more general, if $\Theta$ is closed,
it defines a parallel vectorfield. Of course, this is not always true if
the recurrent vector field is lightlike. Hence, if we use the term `recurrent' we always mean `recurrent
and lightlike'.
A Lorentzian manifold with lightlike parallel vector field is called {\em
Brinkmann-wave}, due to \cite{brinkmann25}.
One has the following description in coordinates (see for example \cite{bb-ike93}).

\blem
$(M,h)$ be a Lorentzian manifold of dimension $n+2$  with recurrent vector field
if and only if there are local  
 coordinates $(U,\varphi=(x, (y_i)_{i=1}^n, z)) $ in which the metric $h$ has the
form
%\begin{equation}\label{walker}
\[h = 2\ dx dz + \sum_{i = 1}^n u_i dy_i dz  + f dz^2 +  \sum_{i,j =
1}^n g_{ij} dy_i\ dy_j,\text{ with $ \frac{\partial g_{ij}}{\partial x}= \frac{\partial u_i}{\partial x}=0$, }   
\]
%\end{equation}
and $f\in
C^\infty(M)$ obeying $\frac{\partial f}{\partial x}=0$ if and only the recurrent vector field can be rescaled to a parallel one. In this case the coordinates can be chosen such that $u_i=0$ and end even that $f=0$ \cite{schimming74}. 
\elem

Before we consider special classes of Lorentzian manifolds with recurrent vector fields we should make some remarks on algebraic properties of the {\em metric} holonomy group. 
The holonomy algebra $\mf{h}$ of a $(n+2)$-dimensional  Lorentzian manifold with recurrent vector field is contained in 
the parabolic algebra $(\rr\+\lason )\ltimes \rrn$. Its projection on $\rrn$ is surjective if and only if the holonomy representation  is {\em indecomposable} (i.e. admits no {\em non-degenerate} invariant sub-space). It is Abelian if and only if it is contained in $\rrn$. The recurrent vector field is parallel if and only if the holonomy is
contained in $\lason\ltimes \rrn$. 
The $\lason$--part of the holonomy is called {\em screen holonomy} because it corresponds to the holonomy of the 
so-called {\em screen bundle} $X^\bot/X\rightarrow M$ \cite{thesis}.
There are four different algebraic types of holonomy algebras (see \cite{bb-ike93}),
two of them {\em uncoupled}, being equal to $\lag\ltimes \rrn$ or $\mf{h}=(\rr\+\lag)\ltimes \rrn$, and two with a coupling
between the center of the screen holonomy and the $\rr$-- resp. the $\rrn$--part.
In \cite{leistner02}, \cite{leistner03}, \cite{leistner03b} we showed that the screen holonomy has to be a Riemannian holonomy algebra,
 a fact which yields a classification of holonomy groups of indecomposable, non-irreducible Lorentzian manifolds.
For Lorentzian manifolds with recurrent vector field $X$ we will use a basis  $(X,E_1, \ldots , E_n,Z)$ with
\beq\label{basis}
h(X,Z)=1,\ h(Z,Z)=0,\ E_i\in X^\bot\text{ with } h(E_i,E_j)=\delta_{ij}\text{ and }h(E_i,Z)=0.\eeq
Its curvature satisfies
$
{\cal R}(X,Y)\ =\ 0\ \text{ for any }Y\in X^\bot
$
because
\[
{\cal R}(X,Y,U,V) \ =\ {\cal R}(U,V,X,Y)
\ =\ 
h(\underbrace{{\cal R}(U,V)X}_{\in \rr\cdot X}  ,Y)
\ =\ 0.
\]
This implies 
\beq\label{null1}
 Ric(X,Y)\ =
\ {\cal R}(X,X,Z,Y)
+
{\cal R}(Z,X,X,Y)
+ \sumi {\cal R}(E_i,X,E_i,Y)
=0
\eeq
for the Ricci tensor and any $Y\in X^\bot$. A semi-Riemannian manifold is called {\em Ricci isotropic} or {\em with totally isotropic Ricci tensor} if the image of the Ricci endomorphism is totally isotropic. 
\bs\label{riciso}
Let $(M,h)$ be a Lorentzian manifold with lightlike, recurrent vector field $X$. $(M,h)$ is Ricci-isotropic if and only if $Ric(Y,.)=0$ for any $Y\in X^\bot$.
In particular, an  isotropic Ricci tensor implies that $S=0$.
\es
\bprf
One direction is trivial: if $Ric(Y,.)=0$ for any $Y\in X^\bot$, then $Ric(U)=Ric(U,Z)\cdot X$, i.e. $h(Ric(U),Ric(V))=0$.
Suppose on the other hand that 
\be
0&=& h(Ric(U),Ric(V))\\
&=& Ric(U,X)\cdot Ric(V,Z)+Ric(V,X)\cdot Ric(U,Z)+\sumi Ric(U,E_i)\cdot Ric(V,E_i)
.\ee
From this equation we get by 
 (\ref{null1}) 
\[
0=
h(Ric(Y),Ric(Y))=
\sumi Ric(Y,E_i)^2
\]
for
$Y\in X^\bot$,
and thus $Ric|_{X^\bot\times X^\bot}=0$.
Furthermore it is
\[0=h(Ric(X),Ric(Z))=Ric(X,Z)^2,\] and
\[0=h(Ric(Z),Ric(Z))=\sumi Ric(Z,E_i)^2,\] hence $Ric(Z,Y)=0$ for $Y\in X^\bot$.
Finally it is $S=2\cdot Ric(X,Z)+\sumi Ric(E_i,E_i)=0$ for an isotropic Ricci tensor.
\eprf
%\bs
%Let $(M,h)$ be a Lorentzian manifold with recurrent vector field $X$ and no parallel vector field. Then there is no conformal scale $\vf$ such that the manifold $(M,e^{2\vf}\cdot h)$ admits a lightlike parallel vector field.
%\es
We do not know whether the existence of a recurrent lightlike vector field and a totally isotropic Ricci tensor implies the existence of a parallel lightlike vector field, but in the next section we will give an example where this implication is true.

\section{pp- and pr-waves and their metric holonomy}

Firstly, we want to recall the conventional definition of a pp-wave. 
\begin{de}\label{ppdef}
A {Brinkmann-wave} is called {\em $pp$-wave} if its curvature tensor ${\cal R}$
satisfies the
 trace condition
%\begin{equation}
$tr_{(3,5)(4,6)} ( {\cal R} \otimes {\cal R} ) =0$.
% \label{ewspur}.
%\end{equation}
\end{de}

R. Schimming  proved the following coordinate description and equivalences.

\blem \cite{schimming74}
A Lorentzian manifold $(M,h)$ of dimension $n+2>2$ is a pp-wave if and only if 
there exist local coordinates
 $(U,\varphi=(x, (y_i)_{i=1}^n, z)) $ in which the metric $h$ has the
form
\begin{equation}
\label{ppform}h = 2\ dx dz   + f dz^2 +  \sum_{i =
1}^n  dy_i^2 \ \mbox{, with $ \frac{\partial f}{\partial x}= 0$.} 
\end{equation}

\elem

\blem \cite{schimming74}\label{eqs}
A Brinkmann wave $(M,h)$ with parallel lightlike vector field $X$ is a pp-wave if and only if
one of the following conditions --- in which $\xi$ denotes the 1-form $h(X,.)$ ---
is satisfied:
\begin{enumerate}
\item \label{eq1}
$ \Lambda_{(1,2,3)} \left(\xi \otimes {\cal R}\right) =0$
\item\label{eq2}
There is a symmetric  $(2,0)$-tensor $r$, with $r(X,.)=0$, such that 
$ {\cal R} =  \ \Lambda_{(1,2)(3,4)}\left( \xi\otimes r \otimes \xi\right)$.
\item\label{eq3}
There is a function $\rho$, such that
$ tr_{(1,5)(4,8)} ( {\cal R} \otimes {\cal R} ) = \rho\ \xi \otimes\xi  \otimes 
\xi \otimes \xi$.
\end{enumerate}
\elem

Now we will give  another equivalence for the definition which seems to be simpler than any of
the trace conditions and which makes a generalisation easier. We denote by
$X^\bot$ the parallel distribution of codimension 1, spanned by tangent vectors orthogonal to the recurrent vector field $X$. $\rr\cdot X$ denotes the distribution spanned by $X$.

\bs
A Brinkmann-wave $(M,h)$ with parallel lightlike vector field $X$
 is a pp-wave if and only if
its curvature tensor satisfies:
\begin{equation}
\label{ppeinfach1}
\cur (U,V): X^\bot \longrightarrow \rr \cdot X \mbox{ for all }U,V\in TM,
\end{equation}
or equivalently 
\begin{equation}
\label{ppeinfach2}
\cur (Y_1,Y_2)=0 \mbox{ for all } Y_1,Y_2\in X^\bot .
\end{equation}
\es

\bprf
The equivalence of (\ref{ppeinfach1}) and (\ref{ppeinfach2}) is obvious:
$\cur (Y_1,Y_2)=0$ for all $Y_1,Y_2\in X^\bot$ $\iff$ $\cur (Y_1,Y_2,U,V)=0$ for all $U,V\in TM $ $\iff$ 
$h(\cur (U,V)Y_1,Y_2)=0$ $\iff$    $\cur (U,V)Y_1\in \rr X$.

To check the defining trace condition for pp-waves we fix a basis $(X=X_p, E_1, \ldots , E_n, Z)$  of the form 
(\ref{basis}). Since $X$ is parallel
it is for $Y_i\in T_p M$: 
\begin{eqnarray*}
%\lefteqn{
tr_{(3,5)(4,6)} ( {\cal R} \otimes {\cal R} )_p(Y_1,Y_2,Y_3,Y_4)
%=}\\
&=& \trace \big(h_p\left({\cal R}_p(Y_1,Y_2)\  .\  ,{\cal R}_p(Y_3,Y_4)\  .\ \big)\right)\\
&=&
%h_p\Big(\underbrace{{\cal R}_p(Y_1,Y_2) X}_{= 0} ,
%{\cal R}_p(Y_3,Y_4) Z\Big)
%+h_p\Big(\underbrace{ {\cal R}_p(Y_3,Y_4)X}_{= 0} ,
%{\cal R}_p(Y_1,Y_2) Z\Big)\\
%&&{ }
%+
\sumk h\big({\cal R}(Y_1,Y_2) E_k ,
{\cal R}(Y_3,Y_4) E_k\big)
\end{eqnarray*}
If on one hand the condition (\ref{ppeinfach1}) is satisfied, then 
$h\big({\cal R}(Y_1,Y_2) E_k ,
{\cal R}(Y_3,Y_4) E_k\big) =0 $ because $E_k\in X^\bot$. Hence the trace condition is satisfied.
On the other hand, suppose that
 the trace vanishes. Since 
$ h\big({\cal R}(Y_1,Y_2) E_k ,
{\cal R}(Y_1,Y_2) E_k\big)\ge 0$ for all $k=1, \ldots n$, this implies 
that $ h\big({\cal R}(Y_1,Y_2) E_k ,
{\cal R}(Y_1,Y_2) E_k\big)= 0$ for all $k=1, \ldots n$. But this is 
(\ref{ppeinfach1}).
\eprf 

From this description we obtain easily the Ricci- and scalar curvature of a pp-wave.

\bfolg
A pp-wave is  Ricci-isotropic and has vanishing scalar curvature.
\efolg

Furthermore, the the following fact, which we haved proved in \cite{leistner01}, becomes obvious.
\bfolg
A Lorentzian manifold with recurrent lightlike vector field has Abelian holonomy
if and only if it is a pp-wave.
\efolg

For sake of completeness we shall mention two subclasses of pp-waves. 
The first are the plane waves which are pp-waves with quasi-recurrent curvature, i.e. 
$\nabla \cur = \xi \otimes \tilde{\cur}$ where $\xi=h(X,.)$ and $\tilde{\cur}$ a
$(4,0)$-tensor.
% with the symmetries of the curvature tensor. 
For plane waves the function $f$ in the 
local form of the metric is of the form $f=\sumij a_{ij} y_i y_j$ where the $a_{ij}$ are functions of $z$.
A subclass of plane waves are the Lorentzian symmetric spaces with solvable transvection group, the so-called Cahen-Wallach spaces (see \cite{cahen-wallach70}, also \cite{bb-ike93}). Here the function $f$ satisfies
 $f=\sumij a_{ij} y_i y_j$ where the $a_{ij}$ are constants.

\bigskip

Now we introduce a new class of non-irreducible Lorentzian manifold by supposing (\ref{ppeinfach1}) but only the existence of a recurrent vector field.
Assuming that the abbreviation `pp' stands for `plane fronted with parallel rays' we shall call them {\em
pr-waves}: `plane fronted with recurrent rays'.
%\footnote{We are not sure if the name 'wave' is still
%justified. {\em Screen-flat} manifolds might be a better name, in accordance to the fact that
%the $\lason$--projection of the holonomy is the holonomy of the {\em screen bundle}.  
%Comments and suggestions for naming these manifolds are very welcome.}

\begin{de}\label{prdef}
We call a Lorentzian manifold $(M,h)$ {\em 
pr-wave} if it admits a recurrent vector field $X$ and   its curvature tensor ${\cal R}$ obeys
\begin{equation}\label{preq}
\cur (U,V): X^\bot \longrightarrow \rr \cdot X \mbox{ for all }U,V\in TM,
\end{equation}
or equivalently $
\cur (Y_1,Y_2)=0 \mbox{ for all } Y_1,Y_2\in X^\bot .$
\end{de}

Since $X$ is not parallel all the trace conditions which were true for a pp-wave,  fail to
hold for a pr-wave. 
 For example, if we suppose (\ref{preq}) we get for the trace
$tr_{(3,5)(4,6)}(\cur \otimes \cur )(U,V,W,Z)=
h_p(\cur (U,V)X, \cur(W,Z)Z)$ which is not necessarily zero.
 But we can prove an equivalence similarly to \ref{eq1} of
Lemma \ref{eqs}.

\blem
A Lorentzian manifold $(M,h)$ with recurrent vector field $X$ is a  
pr-wave if and only if 
$ \Lambda_{(1,2,3)} \left(\xi \otimes {\cal R}\right) =0$, where $\xi$ denotes again the
1-form $h(X,.)$.
\elem
\bprf
Suppose that $(M,h)$ is a pr-wave, fix a vector $Z\in T_p M$ with $h_p(X_p,Z)=1$
 and consider the skew-symmetrisation for $U,V,W\in TM$:
\[\xi(U)\cur (V,W)+\xi(V)\cur (W,U)+\xi(W)\cur (U,V).\]
If $U,V,W\in X^\bot$ this expression is zero of course. For $U,V\in X^\bot$ and $W=Z$ it is equal
to $\xi(Z)\cur (U,V)$, but this is zero because of (\ref{ppeinfach2}).
In case that only $U\in X^\bot$ and $V=W=Z$ it is equal to 
$\xi(Z)\left( \cur (Z,U)+\cur (U,Z)\right)$, which is zero because of the skew-symmetry of the curvature.
On the other hand the vanishing of the skew-symmetrisation implies that 
$\xi(Z)\cur (U,V)=0$ which gives (\ref{ppeinfach2}).
\eprf

Also we get a similar description in terms of local coordinates as for pp-waves.
\blem 
A Lorentzian manifold $(M,h)$ of dimension $n+2>2$ is a pr-wave if and only if 
around any point $o\in M$ exist coordinates
 $(U,\varphi=(x, (y_i)_{i=1}^n, z)) $ in which the metric $h$ has the following
form,
\begin{equation}
\label{prform}h = 2\ dx dz   + f dz^2 +  \sum_{i =
1}^n  dy_i^2  \ \mbox{, with $ f\in C^\infty(M)$.}
\end{equation}

\elem
The {\em proof} is similar to the proof of \cite{schimming74} for pp-waves.
As for pp-waves we can show the relation to the holonomy.

\bs
Let $(M,h)$ be a Lorentzian manifold  with recurrent vector field. The screen holonomy 
$\lag:=pr_{\lason}\left(
\mf{hol}_p(M,h)\right)$ is zero (i.e. $(M.h)$ has solvable holonomy contained in $\rr\ltimes \rrn$) 
if and only if $(M,h)$ is a pr-wave.
\es

\bprf
Both directions follow from the Ambrose-Singer holonomy theorem.
One direction is trivial: If the holonomy of $(M,h)$ is contained in $\rr\ltimes \rrn$, then
--- since $R(U,V)\in\mf{hol}_p(M,h) $  for any $U,V\in T_p M$ --- we get easily the relation
(\ref{ppeinfach1}). 
On the other hand,
let $(M,h)$ be a pr-wave.
Fix a basis $(X,E_1, \ldots ,E_n,Z)$ in $T_p M$ as in (\ref{basis}) and set $E:=
span (E_1, \ldots, E_n)$. Then, by the holonomy theorem,  $\lag:=pr_{\lason}\left(
\mf{hol}_p(M,h)\right)$ is generated  by the following endomorphisms of $E$:
\[pr_E\circ {\cal P} ^{-1}\circ {\cal R}( {\cal P}(U), {\cal P}(V)) \circ {\cal P}_{|E},
\]
$ {\cal P}$ being
the parallel displacement along a piecewise smooth curve starting in $p$, and $U,V\in T_p M$.
Since $X^\bot $ is invariant under parallel displacements, (\ref{ppeinfach1}) ensures that
${\cal P} ^{-1}\circ {\cal R}( {\cal P}(U), {\cal P}(V)) \circ {\cal P}$ maps
$X^\bot$ onto $\rr\cdot X$, hence the $\lason$-part of the holonomy is zero.
\eprf

Finally, we see that Ricci-isotropy forces a pr-wave to be a pp-wave.

\bs
A pr-wave is a pp-wave if and only if it is Ricci-isotropic.
\es
\bprf
We have to show that the recurrent vector field $X$ which satisfies $\nabla X=\Theta \otimes X$ can be rescaled to a parallel one. 
This is possible if $\Theta $ is a differential, i.e. if it is closed. But this is the case if
\be 
0&=&\cur (U,V) X\
=\ d\Theta (U,V)\cdot X
\ee
for $U,V\in TM$. 
Hence we have to show that $\cur (X,U,V,W)=0$ for $U,V,W\in TM$. This is always the case if $U\in X^\bot$ and, for a pr-wave, if $V,W\in X^\bot$. 
But as the Ricci tensor is isotropic it is by proposition \ref{riciso}
\be 
0&=&
Ric(Y,Z)
\\
&=&
\cur(Z,Y,X,Z)+\sumi \underbrace{\cur (E_i,Y,E_i,Z)}_{=0\text{ for pr-waves}}
\ee
where $Y\in TM$ and $(X,E_1, \ldots , E_n,Z)$ a basis of $T_p M$ as in (\ref{basis}).
\eprf

\section{Proof of Theorem \ref{theo}}

Now we turn to the proof of Theorem \ref{theo}. First we prove the easy part.

\bs
Let $(M,c)$ be a simply connected manifold with conformal structure $c$ of Lorentzian signature. If there is a metric in the conformal class $c$ which admits a lightlike recurrent vector field and a totally isotropic Ricci tensor, then the conformal tractor holonomy admits a $2$-dimensional, totally isotropic, invariant subspace.
\es
\bprf
Let $h\in c$ be the Lorentzian metric with recurrent lightlike vector field $X$ and totally isotropic Ricci tensor. Then, by proposition \ref{riciso}, it is  $S=0$, and hence $P=Ric$,  $Ric(Y)=0$ if $Y\in X^\bot$ and $Ric(U)=Ric(U,Z)\cdot X$ else, for $Z$ transversal to $X^\bot$ and $h(X,Z)=1$. Using the decomposition of the tractor bundle with respect to $h$ we consider the following sub-bundle
\[\cal H\ :=\ \cal E [1]\+\cal X,\]
where $\cal X=\rr\cdot X$ is the isotropic line sub-bundle of $TM$ generated by the recurrent vector field $X$. This bundle is totally isotropic and is left invariant by the tractor connection: for $Y\in \Gamma(\cal X)$ and $\s\in \Gamma(\cal E[1])$ we get
\[ D_U (\s, Y, 0)\ =\ \left(U(\s)-h(Y,U), \underbrace{\nabla_U Y + \s  Ric(U,Z)\cdot X}_{\in \rr\cdot X} , -\underbrace{P(Y,U)}_{=0}\right)\in \cal H\ .\]
This proves that $\cal H$ is holonomy invariant. 
\eprf

Now we prove the other direction of theorem \ref{theo}.

\btheo
Let $(M,h)$ be a Lorentzian manifold, ${\cal T}$ its conformal tractor bundle and $D$ the conformal tractor connection. If the holonomy group of the $D$ admits a 2-dimensional totally isotropic invariant subspace, then $(M,h)$ is locally conformally equivalent to a Lorentzian manifold with recurrent lightlike vector field and  totally isotropic Ricci tensor.
\etheo

\bprf
We fix a metric $h$ in the conformal class $c$ inducing the splitting of the tractor bundle
${\cal T}={\cal E}[1]\+TM[-1]\+{\cal E}[-1]$. 
Let ${\cal H}\subset {\cal T} $ be the 2-dimensional, totally isotropic, holonomy-invariant sub-bundle of ${\cal T}$, and ${\cal H}^\bot$ its orthogonal complement with respect to $\la.,.\ra$, which is holonomy-invariant as well.
In order to define a lightlike vector field on $M$ which will be recurrent with respect to a metric from the conformal class of $h$ we consider the intersection of ${\cal H}$ and ${\cal H}^\bot$ with $TM[-1]\+{\cal E}[-1]$, i.e.
\be
\tilde{\cal X}& := & {\cal H}\cap \left(TM[1]\+{\cal E}[-1]\right)\ =\ \{(0,X, \rho)\in {\cal H}\},\\
\tilde{\cal X}^\bot& := & {\cal H}^\bot\cap \left(TM[1]\+{\cal E}[-1]\right)\ =\ \{(0,Y, \rho)\in {\cal H}^\bot\}.
\ee
Since we can assume that ${\cal E}[-1]$ is not contained in ${\cal H}$ (see lemma \ref{clem}), on one hand
$\tilde{\cal X} $ is one-dimensional  and $\tilde{\cal X}^\bot$ is $(n+1)$-dimensional, and on the other hand
the projection ${\cal X}:=pr_{TM} \tilde{\cal X}=\{X\in TM | (0,X,\rho)\in \tilde{\cal X}\} $ is an isotropic line. Its orthogonal complement is the lightlike sub-bundle 
${\cal X}^\bot = pr_{TM} \tilde{\cal X}^\bot$ of codimension $1$ in $TM$.

Of course ${\cal X}^\bot $ is integrable:
if we choose a frame field $(X,E_1, \ldots , E_n) $ of ${\cal X}^\bot$ such that 
\beq\label{cbasis}
X\in \Gamma({\cal X}) \ \text{ and }\ h(E_i, E_j)=\delta_{ij},\eeq then 
 we get \[D_{E_j}(0,E_i,\rho)
=
(0,\nabla_{E_j}E_i+ \rho E_j, ...)
\in \Gamma(\tilde{\cal X}^\bot),\]
for $i\not= j$ and $(0,Y_i,\rho)\in \Gamma(\tilde{\cal X}^\bot)$.
Hence $\nabla_{E_j}E_i\in {\cal X}^\bot$ i.e. $\left[E_i,E_j\right]\in {\cal X}^\bot$. This implies that ${\cal X}^\bot$ is an integrable distribution, the leaves of which being lightlike hypersurfaces.

Since the definition of ${\cal X}$ and its orthogonal ${\cal X}^\bot$ is conformally invariant  we can change the metric in the conformal class in order to perform our calculations.
In order to specify a metric  we have to define the equivalent of a second fundamental form of ${\cal X}^\bot$. Therefore we extend the basis of ${\cal X}^\bot$ to a basis $(X,E_1, \ldots , E_n, Z)$ of $TM$ satisfying  
\beq\label{basis2}
h(X,Z)=1\ ,\ \ h(Z,E_i)=0\ \text{ and  (\ref{cbasis}).}
\eeq
 Then we define the  second fundamental form of ${\cal X}^\bot$ with respect to the transversal vector $Z$ (see for example \cite{bejancu-duggal})
\beq
\begin{array}{rcrcl}
S^Z&:& \cal X^\bot \times \cal X^\bot&\rightarrow & \cal Z:= \rr\cdot Z\\
&&(U,V)&\mapsto &pr_{\cal Z} \nabla_U V= h(\nabla_U V, X)\cdot Z.
\end{array}
\eeq
Obviously it is $S^{f\cdot Z}(U,V)=h(\nabla_UV,\frac{1}{f}X)\cdot f Z= h(\nabla_UV,X)\cdot  Z=S^Z(U,V)$.

Now we claim that there is a metric in the conformal class $c$  and a transversal vector $Z$ defining a basis $(X,E_1, \ldots , E_n, Z)$ such that
\beq\label{claim}\trace S^Z:=\sumi S^Z(E_i,E_i)\ =\ 0.\eeq
To prove this we start with a metric $h$ and fix local coordinates $(x, y_1, \ldots , y_n, z)$ on $M$ such that $x$ parametrizes ${\cal X}$ and $(x, y_1, \ldots , y_n)$ parametrise ${\cal X}^\bot$. We set $Z=\ddz$ and fix a frame field $(X, E_1, \dots , E_n, Z)$ of the form (\ref{basis2}). We define a function $f$ on the coordinate neighborhood by
\[\sumi S^{Z}(E_i,E_i)\ =\ f\cdot Z.\]
Then we define a conformal scale $\vf$ as the solution of the  differential equation
\[d\vf (Z)=\frac{\partial \vf}{\partial z}\ =\ \frac{f}{n},\]
which can be solved by using characteristics. Then we
consider the conformally changed metric $\tilde{h}=e^{2\vf} \cdot h$. 
A basis of the form (\ref{basis2}) is $(e^{-2\vf}X,e^{-\vf}E_1,  \ldots , e^{-\vf}E_n, Z)$.
Its second fundamental form satisfies
\be
\tilde{S}^Z(U,V)&=&
\tilde{h}\left(\widetilde{\nabla}_U V, e^{-2\vf}X\right)\cdot Z\\
&=&
\left( h(\nabla_U V, X) +d\vf (U)\underbrace{h(V,X)}_{=0}
+d\vf (V)\underbrace{h(U,X)}_{=0}
-h(U,V)h(\mathsf{grad}\vf,X)\right)\cdot Z\\
&=&
S^Z(U,V)-h(U,V)\cdot d\vf(Z) Z.
\ee
Then the way $\vf$ was chosen ensures that
\be
\trace  \tilde{S}^Z&=&\sumi \tilde{S}^Z(\tilde{E}_i, \tilde{E}_i)
%\ =\
%e^{-2\vf} \tilde{S}^Z(E_i,E_i)
\\&=&
e^{-2\vf}\left( S^Z(E_i,E_i)-n\cdot d\vf(Z) Z\right)
\\&=&
e^{-2\vf}\big( f-n\cdot d\vf(Z) \big)\cdot Z
\\&=& 0.
\ee
This proves the claim.

From now on we fix the metric in which (\ref{claim}) is satisfied.
Now take 
 $(0,U,\rho)\in \Gamma(\tilde{\cal X})$ and $V\in {\cal X}^\bot$. Then
\[
D_V(0,U,\rho)\ =\ \left(0\ ,\nabla^h_V U+\rho V\ , V(\rho)-P(U,V)\right)
\ \in \ {\cal H}^\bot \cap \left(TM[1]\+{\cal E}[-1]\right) = \tilde{\cal X}^\bot,
\] since ${\cal H}^\bot$ is holonomy-invariant.
In particular it is 
\beq\label{zwei} \nabla^h_V U\in {\cal X}^\bot,
\eeq
which shows that $S^Z(U,V)=0$ for $U\bot V\in {\cal X}^\bot$.
On the other hand we get 
 \[
D_Y(0,X,\rho)\ = \ \left(0\ ,\nabla^h_Y X+\rho Y\ , Y(\rho)-P(X,Y)\right)
\ \in \ {\cal H}\cap \left(TM[1]\+{\cal E}[-1]\right) = \tilde{\cal X},
\]
for 
 $(0,X,\rho)\in \Gamma(\tilde{\cal X})$ and $Y\in {\cal X}^\bot $, because ${\cal H}$ is holonomy-invariant.
In particular it is $\nabla^h_Y X+\rho Y\sim X$ which implies that
\beq\label{eins}
\rho h(Y,Y)\ =\ -h(\nabla^h_Y X,Y)\ =\ -h(\nabla^h_Y Y,X), 
\eeq
since $h(X,Y)=0$.
But it is $h(\nabla^h_Y Y,X)Z=-S^Z(Y,Y)$ and thus
$S^Z(Y,Y)=-\rho\  h(Y,Y) Z$.
By the  choice of the metric above this implies that 
\[0=\sumi S^Z(E_i,E_i)=\rho \cdot n\cdot Z,\]
 yielding $\rho=0$. 
We obtain one one hand
\beq\label{drei}
\begin{array}{rcll}
\nabla^h_Y X & \in &\Gamma( {\cal X})& \mbox{ for }Y\in {\cal X}^\bot \text{ and }X\in \Gamma({\cal X})\ \text{ and }\\
\nabla^h_Y U & \in & \Gamma({\cal X}^\bot)& \mbox{ for }Y\in {\cal X}^\bot \text{ and }U\in \Gamma({\cal X}^\bot),
\end{array}
\eeq
and on the other hand
$
 \{\left(0,X,0\right)|X\in {\cal X}\}=\tilde{\cal X}\subset {\cal H}.
 $

In a next step we show that the covariant derivative in the transversal direction $Z$, which shall be fixed, leaves $\Gamma({\cal X})$ and $\Gamma({\cal X}^\bot)$ invariant.
First we notice that $\Gamma({\cal H})$ contains $(0,X,0)$ and 
$D_Z (0,X,0)=\left(-1, \nabla^h_Z X,-P(X,Z)\right)$ for $X\in \Gamma({\cal X})$. Since ${\cal H} $ is totally isotropic this gives
$0=\left\la(0,X,0),  D_Z (0,X,0)\right\ra=h(X,\nabla^h_Z X)$, i.e. $\nabla^h_Z X\in \Gamma({\cal X}^\bot)$.
Secondly,
 we
get 
\be
0&=&
\left\la D_Z (0,X,0),(0,V,\rho)\right\ra\\
&=&
\left\la (-1, \nabla^h_Z X, -P(X,Z)), (0,V,\rho)\right\ra
\\&=&
-\rho
+h(\nabla^h_Z X,V)
\ee
for $(0,V,\rho)\in\tilde{\cal X}^\bot$ and $(0,X,0)\in \Gamma({\cal X})$, i.e.
\beq\label{vier}
\rho\ =\ 
h(\nabla^h_Z X,V).
\eeq
Considering the second derivative
\[
D_U D_Z (0,X,0)\ =\ \left(-h(\nabla^h_Z X,U)\ ,\  \nabla^h_U \nabla^h_Z X-
P(U)^\sharp-P(X,Z) U\ ,\  ...\right)\in {\cal H},\]
and 
pairing this with $(0,V,\rho)\in {\cal H}^\bot$ we get 
\be
0&=&
\left\la D_U D_Z (0,X,0),(0,V,\rho)\right\ra\\
&=&
-\rho h(\nabla^h_Z X,U) + h(\nabla^h_U\nabla^h_Z X,V) -P(U,V)
\\
&=&
-h(\nabla^h_Z X,V)
h(\nabla^h_Z X,U) + h(\nabla^h_U\nabla^h_Z X,V) -P(U,V).
\ee
Hence the bilinear form
\[(U,V)\longmapsto h(\nabla^h_U\nabla^h_Z X,V)= h(\nabla^h_Z X,V)h(\nabla^h_Z X,U)+P(U,V)\]
is symmetric in $U$ and $V$  from ${\cal X}^\bot$.

Now we want to change the metric conformally in a way that in the new metric $\tilde{h}$ it holds that $\nabla^{\tilde{h}}_Z X\in \Gamma({\cal X})$.
Therefore we consider coordinates $(x, y_1, \ldots , y_n, z)$ on $M$ such that $x$ parametrizes ${\cal X}$ and $(x, y_1, \ldots , y_n)$ parametrise ${\cal X}^\bot$.
Then we consider a conformal  scale  $\vf$ depending only on $(y_1, \ldots y_n)$ and defining the new metric $\tilde{h}=e^{2\vf}\cdot h$.
The new covariant derivative is given by
\[\nabla^{\tilde{h}}_U V=
\nabla^h_U V+d\vf(U) \cdot V + d\vf(V) \cdot U + h(U,V) \cdot \mathsf{grad}\ \vf.\]
This shows that the new derivative in ${\cal X}^\bot$ directions $\nabla^{\tilde{h}}_U$, $U\in {\cal X}^\bot$ still leaves ${\cal X}$ and ${\cal X}^\bot$ invariant.
For $X:=\ddx\in \Gamma({\cal X})$ and a vector $Z$ such that $h(X,Z)=1$ and $h(
\frac{\partial}{\partial y_i},Z)=0$ for all $i=1, \ldots , n$ we get
\be
h( \nabla^{\tilde{h}}_Z X,Y)&=&
h(\nabla^h_Z X,Y) - h(\mathsf{grad}\ \vf, Y)
\\
&=&
\left( h(\nabla^h_Z X,.) - d\vf\right) (Y),
\ee
with  $Y\in \mathsf{span }\left(\frac{\partial}{\partial y_1}, \ldots , \frac{\partial}{\partial y_n}\right)$.
Hence, we can choose the scale $\vf$ in a way that this term vanishes if 
$h(\nabla^h_Z X,.)$ is closed considered as a $1$-form on 
$ \mathsf{span }\left(\frac{\partial}{\partial y_1}, \ldots , \frac{\partial}{\partial y_n}\right)$.
 But it is obviously
 \be
 d(h(\nabla^h_Z X,.) )(U,V)&=&
 U\left(h(\nabla^h_Z X,V)\right)-  V\left(h(\nabla^h_Z X,U)\right) -  h(\nabla^h_Z X,[U,V])
 \\&=&
h(\nabla^h_U\nabla^h_Z X,V)-h(\nabla^h_V\nabla^h_Z X,U)
\\ &=&0,
\ee 
since the last term was proven to be symmetric.
Hence, in the new scale we have that 
$\nabla^{\tilde{h}}_Z$ leaves ${\cal X}$ invariant as well.
But this shows that there is a recurrent, lightlike vector field $X\in \Gamma({\cal X})$ on $(M,\tilde{h})$. But by (\ref{vier}) we get 
\[\tilde{\cal X}^\bot =\{(0,Y,0)|Y\in {\cal X}^\bot\} \subset {\cal H}^\bot.\]
This implies  for $U\in TM$ and $(0,Y,0)\in \Gamma (\tilde{\cal X})^\bot$ that
\be D_U (0,Y,0)& =& \left(-h(U,Y), \nabla^h_U Y, -P(U,Y)\right)\ \in\ \Gamma( {\cal H}^\bot).
\ee
%%Thus if $h(U,Y)=0$, $D_U (0,Y,0)\in {cal X}^\bot$ which implies in this case that $P(U,Y)=Ric(U,Y)=0$. 
Furthermore we observe that  for $(0,X,0) \in\Gamma( {\cal H})$ 
it holds  
\be0&= &\left\la 
D_Z(0,X,0), D_Z(0,X,0)\right\ra
\\
&=&
\left\la
(-1,\nabla^h_Z X,-P(X,Z)), (-1,\nabla^h_Z X,-P(X,Z))\right\ra
\\&=&
2 P(X,Z).
\ee
Hence, $D_Z (0,X,0)=(-1, \nabla_Z X, 0)\in {\cal H}$ and finally \beq\label{fuenf}
0\ =\ \left\la D_Z (0,X,0), D_U (0,Y,0))\right\ra
\\
\ =\ 
P(Y,U)+ \underbrace{h(\nabla^h_Z X, \nabla^h_U Y)}_{=0}
\eeq
for $U\in TM$ and $Y\in \Gamma({\cal X}^\bot)$.
As in the proof of proposition \ref{riciso} this equation  shows 
 that the scalar curvature $S$ of the metric $h$ vanishes. 
%This can be seen by using a basis $(X, E_1, \ldots, E_n , Z)$ such that 
%$Z$ is lightlike and transversal to $\cal X^\bot$,
%$X\in \Gamma(\cal X)$, $E_i\in \Gamma(\cal X^\bot)$ orthonormal to each other and orthogonal to $X$ and $Z$ and $h(X,Z)=1$. 
%Then the $\mathsf{trace} P$ equals to
%\[ \mathsf{trace}\  P\ =\ 
%2P(X,Z)+ \sumi P(E_i,E_i)
%\ =\
%0
%\] because of (\ref{fuenf}).
We conclude that $P=Ric$. Then (\ref{fuenf}) implies that $h$ has a totally isotropic Ricci-tensor.
\eprf

\section{Conformal holonomy of plane waves}
Finally we want to deal with a very simple example of Lorentzian manifolds with recurrent lightlike vector field, the {\em plane waves}. For these spaces the lightlike recurrent vector field is parallel, they are pp-waves, thus  their metric holonomy is $\rrn$, and they are totally Ricci isotropic. In fact, they are even conformally Ricci-flat. This can be seen directly by looking at the transformation formula of the Ricci tensor but we will establish this by calculating the conformal holonomy.

For a plane wave exist coordinates $(x, y_1, \ldots , y_n, z)$ such that the metric has the following form
\beq\label{pw}h= 2dx\ dz + \left(\sumij a_{ij} \ y_i y_j\right)dz^2 + \sumi dy_i^2 ,
\eeq
 where the $a_{ij}$ are functions only of $z$. We set 
 $X=\ddx, \ Y_i=\ddi$ and $Z=\ddz$. 
 The only non-vanishing curvature terms are
 $\cal R \left( Y_i, Z,Z,Y_j\right) \ =\ - a_{ij}$,
 which establishes that
 $Ric  = a\ dz^2$, where $ a=\sumi a_{ii}$.
 We obtain the following result for the conformal holonomy of plane waves.
 \bs
A $(n+2)$-dimensional plane wave $(M,h)$ is conformally Ricci-flat, and if it is indecomposable, its conformal holonomy satisfies
\be
\mf{hol}_p(M,[h])&=&
\rr^{2n+1}\\
&=&
\left\{\left. \left (\begin{array}{cccrr}
0&0&u^t&c&0\\
0&0&v^t&0&-c\\
0&0&0&-v&-u\\
0&0&0&0&0\\
0&0&0&0&0
\end{array}
\right)\right| \, 
u,v\in\rr^n, \ c\in \rr\right\}.\ee
In particular, the tractor connection has two isotropic parallel sections, i.e. there are two Ricci-flat metrics which are locally conformally equivalent to $(M,h)$.
\es
\bprf
We consider the isotropic section $
(\s, \tau \cdot X, 0)$ of the tractor bundle. Since $X:=\ddx$ is parallel 
its tractor derivative is
\[D_U(\s,\tau\cdot X, 0)\ =\ \left(U(\s)-\tau h(U,X)\ ,\  \Big(U(\tau)+\frac{a}{n-2}\ dz(U)\Big)\cdot X\ ,\  0\right).\]
This is zero if $\s$ and $\tau $ depend only on the coordinate $z$ and satisfy the system of ordinary differential equations
\be
\s'&=&\tau\\
\tau'&=& \frac{a}{n-2}\ \s
\ee
where $a=a(z)$ smooth. This system has two independent solutions yielding two parallel isotropic sections of the tractor bundle. Both parallel sections give local scales to Ricci-flat metrics.
Hence 
a plane wave is locally conformally Ricci-flat. By corollary (\ref{richol}) its conformal holonomy is contained in $\mf{hol}_p(M,h)\ltimes \rr^{n+1}=\rr^{2n+1}$. We have to show that it is equal to  $\rr^{2n+1}$. We show that the ambient holonomy contains the holonomy of $g$ and that its projection onto $\rr^{n+1}$ is surjective.

By (\ref{pareqn}) and (\ref{aeqn}) the holonomy of the ambient metric 
$\mf{hol}(\ol{M},\ol{h})$ is generated by homomorphisms which send 
$Y_i$ onto $a_i(1) X+  {\cal P}^{-1}_{\gamma(1)}\circ {\cal R}_{\gamma(1)}(U,V) \circ {\cal P}_{\gamma(1)}(Y_i)$ where $\gamma$ is a curve in $M$, $\cal P_\gamma(1)$ the parallel displacement w.r.t. $h$ and
$a_i$ is determined by 
\begin{eqnarray}
 \nonumber
\dot{a}_i(t)&=& 
-h\Big(\dot{\gamma}(1-t), \underbrace{\cal P_{\gamma|_{[0,1-t]}}\left( \cal P_{\gamma(1)}^{-1}\circ\cal R_{\gamma(1)}(U,V)\circ \cal P_{\gamma(1)}(Y_i)\right)}_{\in \rr\cdot X}\Big).
\end{eqnarray}
Hence, for $\gamma$ a curve with constant $z$-component we obtain  $\dot{a}_i\equiv 0$ and thus $a_i\equiv 0$. This means that the ambient holonomy contains ${\cal P}^{-1}_{\gamma(1)}\circ {\cal R}_{\gamma(1)}(U,V) \circ {\cal P}_{\gamma(1)}$ for curves with constant coordinate $z$. To show that this generates the whole holonomy algebra of $h$ we assume the contrary. Since $Hol(h)$ is supposed to be indecomposable this would imply that there is a vector $Y=\sum_k\eta_kY_k(\gamma(0))$ which is mapped onto zero under all ${\cal P}^{-1}_{\gamma(1)}\circ {\cal R}_{\gamma(1)}(U,V) \circ {\cal P}_{\gamma(1)}$ for curves with constant coordinate $z$. But for these curves the parallel displacement of $Y_k$ is the identity, i.e. $Y=\sumk \eta_k Y_k$ is a vector field which satisfies $\cal R(U,V)Y=0$ for all $U,V\in TM$. This is the integrability condition for the existence of a vector field $\nabla_UY=c\cdot U$ for all $U\in TM$ and a constant $c$. But this is impossible for an indecomposable plane wave, which can be seen by calculating its Christoffel symbols. Hence the holonomy of $h$ is contained in the ambient holonomy.
 
Finally we have to show that the projection of the ambient holonomy on $\rr^{n+1}$ is surjective. Therefore we consider in the expression above geodesics of $(M,h)$ starting in $p$. Then for $Y\in \mathsf{span}(Y_1(p), \ldots , Y_n(p))\subset T_pM$ the corresponding $a(t)$ is given by
\[a(t)\ =\ -h\Big(\dot{\gamma}(0),  \cal P_{\gamma(1)}^{-1}\circ\cal R_{\gamma(1)}(U,V)\circ \cal P_{\gamma(1)}(Y)\Big)\cdot t.
\]
If we assume that this is zero we get that  $\cal R_{\gamma(1)}(U,V)\circ \cal P_{\gamma(1)}(Y)=0$ for all geodesics $\gamma$ starting in $p$ and with $h(\dot{\gamma}(0),X(p))\not=0$. This enables us to extend $Y$ to a vector field $\cal P_{\gamma(1)}(Y)$  which satisfies $\cal R(U,V)Y=0$ for all $U,V\in TM$ on an open and dense submanifold of a normal neighbourhood of $p$, and thus on the whole neighbourhood. This is again a contradiction,  which  implies that the $\rr^{n+1}$ projection of the ambient holonomy is $\rrn$. But from the commutator relations in $\rr^{2n+1}$ we obtain that it is the whole $\rr^{n+1}$.
\eprf
We can illustrate this result in the case where the functions $a_{ij}$ in (\ref{pw}) are analytic, or even constant, the latter being equivalent to the property that $(M,h)$ is a Cahen-Wallach space, i.e. a Lorentzian symmetric space with solvable transvection group \cite{cahen-wallach70}. In these cases the ambient metric 
\[\overline{h}=2d\overline{x}d\overline{z} +\overline{z}^2 \left(  2dx\ dz + \big(\sumij a_{ij} \ y_i y_j\big)dz^2 + \sumi dy_i^2 \right)
\]
is analytic as well and we can calculate its holonomy by higher derivatives of the curvature. Although $\overline{\cur}(\overline{Z},.,.,.)=0$ 
%by (\ref{formelcur}), 
we get for the higher derivatives 
%by (\ref{formelcur}) and (\ref{formelz})
\barr{rcccccl}
(\overline{\nabla}_{Y_i}\overline{\cur})(Y_j,Z,Z,\overline{Z})&=& 
\overline{\cur}(Y_j,Z,Z,\overline{\nabla}_{Y_i}\overline{Z})
&=& 
\overline{\cur}(Y_j,Z,Z, \frac{1}{\overline{z}} Y_i)&=& -\frac{a_{ij}}{z},\text{ and}\\[.2cm]
(\overline{\nabla}_{Z}\overline{\cur})(Y_j,Z,Y_i,\overline{Z})\ &=&
\overline{\cur}(Y_j,Z,Y_i,\overline{\nabla}_{Z}\overline{Z})
&=& 
\overline{\cur}(Y_j,Z,Y_i, \frac{1}{\overline{z}} Z)
& = &\frac{ a_{ij}}{z},
\earr
which generate the $\rr$-component, and the additional $\rrn$-component.

\renewcommand{\baselinestretch}{1}

%\bibliography{ALGBIB,thomas,GEOBIB,SPINBIB,HOLBIB,CONF}
\end{document}